\documentclass[12]{amsart}
\usepackage{amsmath,amsfonts,amssymb,amsthm,enumerate,color}
\usepackage[colorlinks=true]{hyperref}
\makeindex
\begin{document}
\title[Representations of Quantum Groups]{Representations of quantum groups defined over commutative rings II}
\author{Ben L. Cox}
\address{Department of Mathematics \\
The Graduate School of the College of Charleston \\
66 George Street  \\
Charleston SC 29424}
\email{ coxbl@cofc.edu}\urladdr{math.cofc.edu/faculty/cox/}
\author{Thomas J.
Enright}
\address{Department of Mathematics \\
University of
California at San Diego, \\
 La Jolla, CA. 92093}
\email{ tenright@cmath.ucsd.edu}
\subjclass{Primary 17B67, 81R10}
\begin{abstract}  In this article we study the structure of highest weight
modules for quantum groups defined over a commutative ring with particular
emphasis on the structure theory for invariant bilinear forms on these
modules.  
\end{abstract}

\numberwithin{equation}{subsection}

\def\scr{\mathcal}
\def\la{\lambda}
\def\e{\epsilon}
\def\Up{\Upsilon}
\def\vpi{\varpi}
\def\ll{\lambda-\varrho}
\def\a{\alpha}
\def\b{\beta}
\def\i{\iota}
\def\p{\prime}
\def\t{\theta}
\def\T{\Theta}
\def\vT{{^\varrho\Theta}}
\def\vPsi{{^\varrho\Psi}}
\def\tiep{T_{e}'}
\def\tiepp{T_{e}''}
\def\twep{T_{w,e}'}
\def\tiempp{T_{-e}''}
\def\J{\mathcal J}
\def\J{\mathcal J}
\def\brho{\varrho}      
\def\ru{{_R\mathbf U}}    
\def\rui{{_R{\mathbf U}}}
\def\bu{\mathbf U}
\def\bui{{{\mathbf U}}}
\def\bz{{\mathbf Z}}
\def\ii{\pi \circ \iota}
\def\ra{\rightarrow}
\def\Par{\mathbb P}
\def\di{\diamond}
\def\vj{w_{\l,j}}
\def\vn{w_{\l,n}}
\def\v(-n){w_{\l,-n}}
\def\vsj{w_{-\l,j}}
\def\vsn{w_{-\l,n}}
\def\v\e(n){w_{\e \l,n}}
\def\v\e(-n){w_{\e \l,-n}}
\def\vs(-n){w_{-\l,-n}}  
\def\U{\mathcal U}           
\def\U_Y{\mathcal U_Y}
\def\gg{\mathfrak g_R}
\def\bb{\mathfrak b_R}
\def\hh{\mathfrak h_R}
\def\pp{\mathfrak p_R}
\def\s{\sigma}
\def\l{\lambda}
\def\r{\rho}
\def\G{\Gamma}
\def\d{\delta}


\newcommand{\qbinom}[2]{\genfrac{[}{]}{0pt}{}{#1}{#2}}

\def\zneg{{\mathbb Z}_{<0}}
\def\znpos{{\mathbb Z}_{\leq 0}}
\def\zpos{{\mathbb Z}_{>0}}
\def\posz{\zpos}
\def\znneg{{\mathbb Z}_{\geq 0}}
\def\vac{\text {\rm Vac}\,}
\def\m#1#2{M_{#1}(#2)}
\def\im{\text {\rm im}\,}
\def\rank{\text {\rm rank}\,}
\def\Tor{\text {\rm Tor}\,}
\def\Hom{\text {\rm Hom}\,}
\def\Span{\text {\rm Span}\,}
\def\Ann{\text {\rm Ann}\,}
\def\End{\text {\rm End}\,}
\def\rad{\text {\rm rad}\,}
\def\inv{\text {\rm Inv}\,}
\def\ind{\text {\rm Ind}\,}

\def\fk #1{{\mathfrak #1}}
\def\hfk #1{\hat {\fk#1}}
\def\dfk #1{\dot {\fk#1}}
\def\angles #1{\langle #1\rangle}
\def\nneg{\fk n_{-}}
\def\npos{\fk n_{+}}
\def\hnpos{{\hat {\fk n}}_{+}}
\def\hnneg{{\hat {\fk n}}_{-}}
\def\bbc{{\mathbb C}}
\def\bbz{{\mathbb Z}}
\def\bbzn{{\mathbb Z}_{<0}}
\def\bbzpos{{\mathbb Z}_{>0}}
\def\bbznpos{{\mathbb Z}_{\leq 0}}
\def\bbznneg{{\mathbb Z}_{\geq 0}}
\def\plusminus{\pm 1}
\def\mla{M(\lambda)}
\def\nla{N(\la)}
\def\bmla{{\overline {M(\lambda)}}}
\def\mhla{M_{\hfk g}(\hat\lambda)}
\def\mmla{M_{\fk m}(\la)}
\def\mkla{M_{\fk k}(\la)}

\def\TT{T^{-1}\cdot}
\def\A{\mathbb A}
\def\B{\mathbb B}
\def\P{\mathbb P}
\def\L{\mathbb L}
\def\D{\mathbb D}
\def\oA{\overline{\mathbb A}}
\def\oB{\overline{\mathbb B}}
\def\oP{\overline{\mathbb P}}
\def\oL{\overline{\mathbb L}}
\def\oD{\overline{\mathbb D}}
\def\hra{\hookrightarrow}
\swapnumbers
\newtheorem{thm}{Theorem}[section]
\newtheorem{lem}[thm]{Lemma}
\newtheorem{cor}[thm]{Corollary}
\newtheorem{prop}[thm]{Proposition}
\newtheorem{rem}[thm]{Remark}
\theoremstyle{definition}
\newtheorem{definition}[thm]{Definition}
\newtheorem{example}[thm]{Example}
\newtheorem{xca}[thm]{Exercise}

\newcommand{\thmref}[1]{Theorem~\ref{#1}}
\newcommand{\propref}[1]{Proposition~\ref{#1}}
\newcommand{\corref}[1]{Corollary~\ref{#1}}
\newcommand{\secref}[1]{\S\ref{#1}}
\newcommand{\lemref}[1]{Lemma~\ref{#1}}
\newcommand{\eqnref}[1]{~(\ref{#1})}
\newcommand{\alignref}[1]{~(\ref{#1})}
\newcommand{\egref}[1]{example~\ref{#1}}
\newcommand{\gatherref}[1]{gather~\ref{#1}}

\dedicatory{Dedicated to V. S. Varadarajan}

\maketitle

\date{March 2001}
\section{Introduction and Summary of Results.} 
Let $v$ be an indeterminate and $k$ a field of characteristic zero. Let
$U$ be the quantized enveloping algebra defined over $k(v)$ with generators
$K^{\pm 1},E,F$ and relations 
\begin{gather*}
[E,F]=\frac{K-K^{-1}}{ v-v^{-1}} , \quad KEK^{-1}=v^2E,\\
 KFK^{-1}=v^{-2}F\quad\text{and} \quad KK^{-1}=K^{-1}K=1.
\end{gather*}
Let $\bu^0$ be the subalgebra generated by $K^{\pm 1}$ and let $B$ be
the subalgebra generated by $\bu^0$ and $E$. More precisely we are
following the notation given in \cite{MR96c:17022} where we take
$I=\{i\}$, $i\cdot i=2$, $Y=\mathbb Z[I]\cong \mathbb Z$, 
$X=\hom(\mathbb Z[I], \mathbb Z)\cong \mathbb Z$, $F=F_i$, $E=E_i$, and
$K=K_i$.

Let $R$ be the power series ring in $T-1$ with coefficients in $k(v)$ i.e. 
$$
R=k(v)[[T-1]]:=\lim_\leftarrow \frac{k(v)[T,T^{-1}]}{ (T-1)^i}.
$$
Set $\mathcal K$ equal to the field of fractions of $R$. Let $s$ be the
involution of $R$ induced by $T\ra T^{-1}$ i.e. the involution that sends
$T$ to $T^{-1}=1/(1+(T-1))=\sum_{i\geq 0}(-1)^i(T-1)^i$. Let the subscript
$R$ denote the extension of scalars from $k(v)$ to $R$ , e.g. $\ru =
R\otimes_{k(v)} \bu$. For any representation $(\pi,A)$ of $\ru$ we can
twist the representation in two ways by composing with automorphisms of
$\ru$. The first is $\pi \circ (s\otimes 1)$ while the second is $\pi
\circ (1\otimes \T)$ for any automorphism $\Theta$ of $\bu$. We designate the corresponding
$\ru$-modules by $A^s$ and $A^{\Theta}$. Twisting the action by both $s$
and $\Theta$ we obtain the composite $(A^s)^{\Theta}=(A^{\Theta})^s$ which
we denote by
$A^{s\Theta}$.

Let $m$ denote the homomorphism of $\ru^0$ onto $R$ with $m(K)=T$. For
$\lambda \in \mathbb Z$ let $m+\lambda$ denote the homomorphism of
$\bu^0$ to
$R$ with $(m+\lambda )(K)=Tv^\lambda$. We use the additive notation
$m+\lambda$ to indicate that this map originated in the classical setting
from an addition of two algebra homomorphisms. It however is not a sum of
two homomorphisms but rather a product. Let
$R_{m+\lambda}$ be the corresponding $_RB$-module and define the Verma
module
 \begin{equation}
_RM(m+\lambda )=\ru\otimes_{_RB} R_{m+\lambda-\rho} .
\end{equation}

Let $\rho_1:\mathbf U\to \mathbf U$ be the algebra isomorphism
determined by the assignment 
\begin{equation}\label{rhoone}
\rho_1(E)=-vF,\quad \rho_1(F)=-v^{-1}E,\quad
\rho_1(K)=K^{-1}. 
\end{equation}  
Define also an algebra
anti-automorphism $\brho:\mathbf U\to \mathbf U$ by 
\begin{equation}
\label{antiauto}
\brho(E)=vKF,\quad \brho(F)=vK^{-1}E,\quad
\brho(K_\mu)=K_\mu. 
\end{equation}
 These maps are related through the antipode $S$ of  $\mathbf U$ by
$\brho=\brho_1S$. 

For $\ru$-modules $M,N$ and $\mathcal F$, let $\mathbb P(M,N)$ and $\mathbb
P(M,N,\mathcal F)$ denote the space of $R$-bilinear maps of $M\times N$ to $R$
and $\mathcal F$ respectively, with the following invariance condition: 
\begin{equation}\label{correction}
\sum x_{(1)}*\phi(Sx_{(3)}\cdot a,\varrho(x_{(2)})b)=\mathbf e(x) 
\phi(a,b)
\end{equation} where $\Delta\otimes 1\circ \Delta (x)=\sum x_{(1)}\otimes
x_{(2)}\otimes x_{(3)}$,  $\mathbf e:\bu\to k(v)$ is the counit and $*$
denote the action twisted by $\rho_1$; in other words $x*n:=\rho_1(x)n$.  If we let
$\hom_\ru(A,B)$ denote the set of module $\ru$-module homomorphisms, then one
can check on generators of $\ru$ that
$\mathbb P(M,N,\mathcal F)\cong \hom_\ru (M\otimes_R N^{\rho_1},{_R\mathcal
F}^\brho)$ (see \cite[3.10.6]{MR96m:17029}).  Formula \eqnref{correction} corrects an
error in  \cite[6.2.2]{MR96c:17022}. Let $\mathbb P(N)=\mathbb
P(N,N)$ denote the
$R$-module of invariant forms on $N$. 

For the rest of the introduction we let $M$ denote the $\ru$ Verma
module with highest weight $Tv^{-1}$ ; i.e. $M=M(m)$ and let $\mathcal
F$ be any finite dimensional $\bu$-module.  A natural parameterization for
$\mathbb P(M\otimes {_R\mathcal F})$ was given in \cite{MR96c:17022}.
Fix an invariant form $\phi_M$ on $M$ normalized as in  ( .). For each
$\ru$-module homomorphism $\b:{_R\mathcal E}\otimes_R\mathcal
F^{\rho_1}\ra\ru$ define what we call the {\it induced form}
$\chi_{\b,\phi_M}$ by the formula, for
$e\in_ R\mathcal E,f\in _R\mathcal F, \ m,n\in N$,
\begin{equation}\label{inducedform}
\chi_{\b,\phi_M}(m\otimes e,n\otimes f)=\phi_M( m,\beta(e\otimes f)*n).
\end{equation}

\begin{prop}
Suppose $\b:{_R\mathcal F}\otimes {_R\mathcal F^{\rho_1}}\ra \ru$ is a module
homomorphism with $\ru$ having the adjoint action.    Then $M\otimes
{_R\mathcal F}$ decomposes as the $\chi_{\b,\phi_M}$-orthogonal sum of
indecomposable $\ru$-modules.
\end{prop}

We now begin the description of our main result:
Recall from \cite{MR96c:17022} we define a {\it cycle} ( for $A$ ) to be a pair $(A,\Psi)$
where
$A$ is a $\mathbf U$ (or $_R\mathbf U$)  module and $\Psi$ is a module
homomorphism  
\begin{equation}
\Psi : A_\pi^{sT_{-1}'} \
\ra \ A .
\end{equation} 
Here $A_\pi$ is defined to be $A_F/\iota A$ where $\iota:A\to A_F$ is the canonical embedding of the module $A$ into its localization $A_F$ with respect to root vector $F$.  Note that modules of the form $A_\pi^{sT_{-1}'}$ above appear naturally in other mathematical work (see \cite{MR2032059} and \cite{MR2074588}) besides our own (see \cite{MR96c:17022}).
We  choose a homomorphism $\Psi : M_F\ra M$ and set
\begin{equation}
\bar \Psi:=\Psi \otimes sT_{1}''\circ L^{-1}:(M_\pi\otimes \mathcal E)^{sT_{-1}'}\to
M\otimes \mathcal E.
\end{equation}
(The linear map $L$ is defined in Lusztig's book - see also \eqnref{L}.)
Let $\iota:\mathbb P(M\otimes \mathcal E, N\otimes
\mathcal F)\to\hom_\ru(M\otimes \mathcal E\otimes (N\otimes
\mathcal F)^{\rho_1},R)$ be the canonical isomorphism with $\iota(\chi)(a\otimes b)=\chi(a,b)$.  Note that $a\otimes b\in
M\otimes
\mathcal E\otimes (N\otimes
\mathcal F)^{\rho_1}$ on the left hand side, while $(a,b)\in M\otimes \mathcal E\times N\otimes
\mathcal F$ on the right hand side.  Define $\chi\mapsto \chi^\sharp$ in $
\End(\mathbb P(M\otimes
\mathcal E, N\otimes
\mathcal F))$ by
\begin{equation}
\iota(\chi^\sharp)(\bar\Psi(a)\otimes \bar\Psi(b)):=
s\circ\iota(\chi_\pi)\circ L(a\otimes b)
\end{equation}
for $a\in(M\otimes \mathcal E)_\pi^{sT_{-1}'}$,
$b\in (N\otimes\mathcal  F)^{s T_{-1}'\rho_1}_\pi$ 
and $\chi\in\mathbb P(M\otimes \mathcal E,N\otimes \mathcal F)$. 

Let  $\mathcal F_m$
and
$\mathcal F_n$ be $X$-admissible finite dimensional $\bu$-modules
given in \secref{CG} with basis $u^{(m)}_k$, $0\leq k\leq m$.  Fix a homomorphism
$\beta:_R\mathcal F_m\otimes _R\mathcal F_n^{\rho_1}\to {_RF(\bu)}$ which has the
form \eqnref{factorization}
\begin{equation}
\beta=\sum_{m,n,k}r^{m,n}_{k}\beta^{m,n}_{2k}
\end{equation}
where $r^{m,n}_{k}\in R$, and $\beta^{m,n}_{2r}$ is defined by 
\begin{equation}
\beta^{m,n}_{2r}(u^{(m+n-2q)})=\delta_{2r,m+n-2q}E^{(r)}
K^{-r}.
\end{equation}
Our main symmetry result on induced forms is \thmref{firstinvariance}:
\begin{thm}
Let $M$ be the Verma module of highest weight $Tv^{-1}$ (so that $\lambda=0$) and assume that $\beta:_R\mathcal F_m\otimes _R\mathcal F_n^{\rho_1}\to
{_RF(\bu)}$  has the form \eqnref{factorization}.
If
$\phi$ is a $\ru$-invariant pairing on $M$ satisfying $s\circ
\phi_{\pi}\circ L=\phi\circ (\Psi\otimes
\Psi)$, then 
\begin{equation}
\chi _{\b,\phi}^\sharp = \chi_{s\beta ,\phi}\ .
\end{equation} 
\end{thm}

Most of the results in sections 1-7 are used in the proof of this theorem. In sections 8 and 9 we give a taste of how one can use induced forms to get information on filtrations of modules.  We plan to pursue this in future work.  

Let $\Pi=\{\alpha,\beta\}$ be the set of simple roots and $\gamma\in\Pi$ for $ \mathfrak g=\mathfrak{sl}(3)$ or $\mathfrak{sp}(4)$.  In the last section we give examples of how one can relate the Shapovalov form for $U_v(\mathfrak g)$, to the Shapovalov form on a reductive subalgebra $U(\mathfrak a)$ generated $E_\alpha,F_\alpha,K_\gamma$, $\gamma\in\Pi$.  In particular we explicitly describe the coefficients $r^{m,n}_{k}$ for particular $\mathbf\beta$ that appear in the study of these Shapovalov forms.  We will expand on this study in future work.

\section{$q$-Calculus.}
\subsection{Definitions}
As many before us have done, for $m\in\mathbb Z$ we define
\begin{align*} 
[m]&:=\frac{v^m-v^{-m}}{v-v^{-1}},  \\
[m]_{(n)}&:=[m]\cdot[m-1]\cdots [m-n+1] \\
[m]!&:=[m]_{(m)}, \quad [0]!:=1  \\
\qbinom{m}{ n}&=\begin{cases}\frac{[m]_{(n)}}{[n]!}&\quad \text{for} \quad n\geq 0 \\ 
  0 &\quad {\rm if}\quad  n<0 . 
\end{cases}
\end{align*}
For $j\geq 0$, Gauss' versions of the Binomial Theorem are 
\begin{equation}\label{gaussbinomial}
\prod_{l=1}^{j}(1-zv^{2(l-1)})=\sum_{k=0}^{j}(-1)^k\qbinom{j}{ k}v^{k(j-1)}z^k
\end{equation}
and 
\begin{equation}\label{gaussbinomial2}
\prod_{l=1}^{j}(1-zv^{2(l-1)})^{-1}=\sum_{k=0}^{\infty}(-1)^k\qbinom{-j}{ k}v^{k(j-1)}z^k
\end{equation}
See \cite[1.157, 1.158]{MR96m:17032}.
For $r\in \mathbb Z$ define
\begin{gather}\label{quantumbinom}
[T;r]:=\frac{v^{r}T- v^{-r}T^{-1}}{ v-v^{-1}},   \\ 
[T;r]_{(j)}:=[T;r][T;r-1]\cdots [T;r-j+1],
\quad[T;r]_{(0)}:= [T;r]^{(0)}:=1,
\quad\text{if}\quad j>0,\notag\\ 
[T;r]^{(j)}:=[T;r+1]\cdots [T;r+j]
\notag \\ \notag \\
\qbinom{T;r}{ j}:=\begin{cases}[T;r]_{(j)}/[j]!
& \text{if}\quad j\geq 0 \\
0&  \text{if}\quad j<0.\end{cases}
\end{gather}
%
Note that $[T;\lambda]^{(k)}$ is invertible in $R$ for $\lambda\geq
0$ and $[T;\lambda+1]_{(k)}$ is invertible provided
$\lambda+1\geq  k$ or $\lambda<0$ ($k\geq 0$). On the other hand $[T;\lambda+1]_{(k)}$ is divisible by $T-1$ and not by $(T-1)^2$ for $k>\lambda+1$ due to
the fact that $[T,r]$ is divisible by exactly $(T-1)^i$ where $i=0$ if $r\neq 0$ and $i=1$ for $r= 0$.  A useful related computation is 
\begin{equation}\label{matrixentry}
[r]![T;r]_{(r)}^{-1}
\equiv \left(1+\left(r-\frac{2}{\{1\}}
\sum_{k=1}^r\frac{v^{k}}{[k]}\right)(T-1)\right)\mod (T-1)^2
\end{equation}
where $\{1\}=v-v^{-1}$.  
Indeed 
\begin{align*}
[T;k]^{-1}
&=T[k]^{-1}\left(
\sum_{i\geq 0}\left(\frac{v^{k}}{1-v^k}\right)^i(T-1)^i\right)\left(
\sum_{j\geq 0}(-1)^j\left(\frac{v^{k}}{1+v^k}\right)^j (T-1)^j\right)
\end{align*}
This implies \eqnref{matrixentry} above.

\subsection{Identities}
Two useful formulae for us will be

\begin{equation}\label{Ma1}
\qbinom{s-u}{ r}=\sum_p
(-1)^pv^{\pm (p(s-u-r+1)+ru)}\qbinom{u}{ p}\qbinom{s-p}{
r-p}
\end{equation}
\begin{equation}\label{Ma2}
\qbinom{u+v+r-1}{ r}=\sum_p
v^{\pm (p(u+v)-ru)}\qbinom{u+p-1}{ p}\qbinom{v+r-p-1}{r-p}
\end{equation}
which come from \cite[1.160a, 1.161a ]{MR96m:17032}, respectively.

We have a yet another variant of the Binomial Theorem:

\begin{lem}\label{thirdbinomial}
\begin{align*}
 v^{k(n+1-k)} [T;2k-n-1]_{(k)}&=\sum_{j=0}^k (-1)^jv^{j(n-2k+1)}
                T^{-j}\qbinom{k}{j} [n-k+j]_{(j)} [T;k]_{(k-j)}\ 
\end{align*}
\end{lem}
\begin{proof}  

The  proof follows from an application of Gauss' binomial theorem and \eqnref{Ma1}
\end{proof}

\begin{lem}[Chu-Vandermonde formula]\label{chuvandermonde}  For integers $k$ and $r$ with $0\leq k\leq r$ we have
\begin{equation}
\sum_{l=0}^{k}(-1)^{l}
  v^{  l (r -k+1) }
    \left[\begin{matrix}T;k \\ l \end{matrix}\right]
    \left[\begin{matrix}T;r+k-l \\ k-l \end{matrix}\right]=v^{-k^2}T^{-k}\qbinom{r}{k}
\end{equation}
\end{lem}

\begin{proof}  
The proof is obtained by using the Taylor series expansion in $T$ together with Gauss' binomial theorem.
 \end{proof}

\section{$\bui$-algebra Automorphisms and Intertwining Maps.} 

\subsection{} Following Lusztig, \cite[Chapter 5]{MR94m:17016}, we let
$\mathcal C'$ denote the category whose objects are $\mathbb Z$- graded
$\bu$-modules $M=\oplus _{n\in\mathbb Z}M^n$ such that 

\begin{enumerate}[(i)]
\item  $E,F$ act locally nilpotently on $M$,
\item  $Km=v^nm$ for all $m\in M^n$.
\end{enumerate}

Fix $e=\pm 1$ and let $M\in \mathcal C'$.  Define Lusztig's automorphisms
$\tiep,\tiepp:M\to M$ by
\begin{equation}\label{defofT}
\tiep(m):=\sum_{a,b,c;a-b+c=n}(-1)^bv^{e(-ac+b)}F^{(a)}E^{(b)}F^{(c)}m, 
\end{equation}
 and

\begin{equation}
\tiepp(m):=\sum_{a,b,c;-a+b-c=n}(-1)^bv^{e(-ac+b)}E^{(a)}F^{(b)}E^{(c)}m
\end{equation}
for $m\in M^n$.  In the above $E^{(a)}:=E^a/[a]!$ is the $a$th {\it
divided power of $E$}.

Lusztig defined automorphisms $\tiep$ and $\tiepp$ on $\bu$ by
$$
T'_e(E^{(p)})=(-1)^pv^{ep(p-1)}K^{ep}F^{(p)},\quad \quad
T'_e(F^{(p)})=(-1)^pv^{-ep(p-1)}E^{(p)}K^{-ep}
$$
and
$$
T''_{-e}(E^{(p)})=(-1)^pv^{ep(p-1)}F^{(p)}K^{-ep},\quad \quad
T''_{-e}(F^{(p)})=(-1)^pv^{-ep(p-1)}K^{ep}E^{(p)}.
$$

One can check on generators that 
\begin{equation}\label{Tandrho}
\rho_1\circ T_{-1}'=T_{-1}'\circ
\rho_1.
\end{equation}
  
If
$M$ is in
$\mathcal C'$,
$x\in \bu$ and $m\in M$, then we have

\begin{equation}
\Theta (x\cdot m) =\Theta(x) \Theta m 
\end{equation}
for $\Theta=\tiep$ or $\Theta=\tiepp$ (see \cite[37.1.2]{MR94m:17016}). The last identity can be interpreted to say that
$\Theta$ and
$\Theta\otimes s$ are intertwining maps;
\begin{equation}
\Theta:M \ra M^\Theta\quad \quad \quad \quad
\Theta\otimes s\ :{_R}M\ra {_R}M^{\Theta\otimes s}.
\end{equation}
To simplify notation we shall sometimes write $s\Theta$
in place of $\Theta\otimes s$. \smallskip

We now describe the explicit action of $\Theta$ on $M$.
\begin{lem}[{\cite[Prop. 5.2.2]{MR94m:17016}}]\label{symmetries} Let
$m\geq 0$ and
$j,h\in [0,m]$ be such that $j+h=m$.
\begin{enumerate}[(a)]
\item If $\eta\in M^m$ is such that $E\eta=0$, then
$\tiep(F^{(j)}\eta)=(-1)^jv^{e(jh+j)}F^{(h)}\eta$.
\item If $\zeta\in M^{-m}$ is such that $F\zeta=0$, then
$\tiepp(E^{(j)}\zeta) =(-1)^jv^{e(jh+j)}E^{(h)}\zeta$.
\end{enumerate}
\end{lem}

Let $F(\bu)$ denote the ad-locally finite submodule of $\bu$.  We know
from \cite{MR95e:17017} that $F(\bu)$ is tensor
product of harmonic elements $\mathcal H$ and the center $Z(U)$.  Here
$\mathcal H=\oplus_{m\in\mathbb Z} \mathcal H_{2m}$ and $\mathcal
H_{2m}=\text{ad}\, \bu (EK^{-1})$.

 There is another category that
we will need and it is defined as follows: Let
$M$ be a
$\ru$-module.  One says that
$M$ is {\it
$_R\bu^0$-semisimple} if $M$ is the direct sum of $R$-modules
$M^{\mu}$ where $K$ acts by $T \ v^\mu$, $\mu\in\mathbb Z$; i.e. by
weight $m+\mu$.  Then $\mathcal C_{R}$ denotes the category of $\ru$-modules
$M$ for which $F$ acts locally nilpotently and  $M$ is
$_R\bu^0$-semisimple.

For $M$ and $N$ two objects in $\mathcal C'$ or one of them is in 
$_R\mathcal C$,
Lusztig defined
the linear map
$L:M\otimes N\to M\otimes N$ given by
\begin{equation}\label{L}
L(x\otimes y)=\sum_n(-1)^nv^{-n(n-1)/2}\{n\}F^{(n)}x\otimes E^{(n)}y
\end{equation} 
%
%
where $\{n\}:=\prod_{a=1}^n(v^a-v^{-a})$ and $\{0\}:=1$.  %
One can show
$$ 
L^{-1}(x\otimes y)=\sum_nv^{n(n-1)/2}\{n\}F^{(n)}x\otimes E^{(n)}y.
$$ 
%

\begin{lem}[\cite{MR94m:17016}]\label{LLemma1}
Let $M$ and $N$ be two objects in $\mathcal C'$.  Then
$ T_{1}'' L(z)=(T_{1}''\otimes T_{1}'')(z)$ for all $z\in
M\otimes N$. 
\end{lem}

\begin{lem}\label{LLemma2}
Let $M$ be a module in $\mathcal C_R$ and $N$ a module in
$\mathcal C'$.  Then for $x\in M^t$ and $y\in N^s$ we have
\begin{align*}
FL(x\otimes y) &=L(x\otimes Fy+v^sFx\otimes y) \\
EL(x\otimes y) &=L(Ex\otimes y+v^{-t}T^{-1}x\otimes Ey)
\end{align*}
\end{lem}

\begin{proof}  Let $\bar{}: k(v)\to k(v)$ be the $\mathbb Q$-algebra
isomorphism given by sending $v$ to $v^{-1}$. More over let
$\bar {}:\bu\to
\bu$ denote the unique
$\mathbb Q$-algebra homomorphism defined by 
$$
\bar E=E, \quad \bar F=F,\quad \bar K_\mu=K_{-\mu},\quad \overline{fu}=\bar f\bar
u,\quad f\in k(v),\,u\in\bu.
$$
From
\cite[Theorem 4.1.2]{MR94m:17016}, we have
$\Delta(u)L=L\overline{\Delta(\bar u)}$ for any $u\in \bu$.  Let
$x\in M^t$ and
$y\in
\mathcal E^s$.  In particular we deduce
\begin{align*}
FL(x\otimes y) &=L\overline{\Delta(F)}(x\otimes y)=L(v^sFx\otimes y+x\otimes Fy).
\end{align*}
Similarly
$EL(x\otimes y) 
  =L (Ex\otimes y+v^{-t}T^{-1}x\otimes Ey)$.
\end{proof}

\begin{cor}\label{Lisomorphism}
Let $M$ be a module in $\mathcal C_R$ and $\mathcal E$ a module in
$\mathcal C'$.  Then $L$ defines an isomorphism of
the $\bu$-module
$M^{T_{-1}'}\otimes \mathcal E^{T_{-1}'}$ onto $(M\otimes\mathcal
E)^{T_{-1}'}$.
\end{cor}

\begin{proof}  The proof follows from a direct calculation on weight vectors in $M$ and $\mathcal E$.
\end{proof} 

\noindent{\bf Remark:}  It seems that one should be able to prove
\corref{Lisomorphism} through the use of \lemref{LLemma1}, however one
must take into account that
$T_e'$ may not be defined on $M$. 

Set 
\begin{equation}\label{mathcalLinv}\index{mathcalLinv}
\mathcal L^{-1}=\sum_p(-1)^pv^{3\frac{p(p-1)}{
2}}\{p\}E^{(p)}K^pF^{(p)},
\end{equation}
and
\begin{align}\label{mathcalL}
\mathcal L
&=\sum_p
    v^{-3\frac{p(p-1)}{2}}\{p\}E^{(p)}K^{-p}F^{(p)}.
\end{align}
and note that $\mathcal L$ and $\mathcal L^{-1}$ are well
defined operators on lowest weight modules.  These maps are similar to some operators defined by Kashiwara.  See for example the 
$\Pi_t$ in \cite[16.1]{MR94m:17016}.

\begin{lem} Suppose that $M$ and $N$ be highest weight
modules with $\psi_M:M^{sT'_{-1}}_\pi\to M$, $\psi_N:N^{sT'_{-1}}_\pi\to
N$ are $\ru$-module homomorphisms and 
$\phi$ a $\rho$-invariant form on $M\times N$.  Then
\begin{equation}\label{Ll}
\phi\circ(\psi_M\otimes \psi_N)\circ L^{-1}=\phi\circ
(\psi_M\otimes \psi_N) \circ (\mathcal  L^{-1}\otimes 1).
\end{equation}

\end{lem}

\begin{proof} Let $m\in M_{\pi}$, $n\in N_{\pi}^{\rho_1}$ and let $*$
denote the action on  $N_{\pi}^{\rho_1}$ so that $x*n=\rho_1(x)n$ where
juxtaposition means the action on $N_\pi$. Recall
$L^{-1}$ is defined on $M_\pi\otimes_R N_\pi^{\rho_1}$, but one can view
it as an $\ru$-module homomorphism from $ (M_\pi\otimes_R
N_\pi^{\rho_1})^{sT_{-1}'}$ to $ M_\pi^{sT_{-1}'}\otimes_R
(N_\pi^{sT_{-1}'})^{\rho_1}$.  Hence the left hand side of \eqnref{Ll}
is in $\Hom_\ru((M_\pi\otimes_R N_\pi^{\rho_1})^{sT'_{-1}},R)$.
Now
\begin{align*}
\phi\circ(\psi_M\otimes \psi_N)\circ L^{-1}&(m\otimes
n)=\sum_{p\geq 0}v^{\frac{p(p-1)}{2}}\{p\}
               \phi(\psi_M(F^{(p)}m), \psi_N(E^{(p)}*n)) \\
&=\sum_{p\geq 0}v^{\frac{p(p-1)}{2}}\{p\}
               \phi(\rho(T''_1\rho_1
                E^{(p)})\psi_N(F^{(p)}m),\psi(n))\\
&=\sum_{p\geq 0}v^{\frac{p(p-1)}{2}}\{p\}
               \phi(\psi_M(T'_{-1}\rho T''_1\rho_1
                (E^{(p)})F^{(p)}m),\psi_N(n))\\
&=\phi\left(\psi_M\left(
             \mathcal  L^{-1}m\right),\psi_N(n)\right)\\
\end{align*}
since
$T'_{-1}\rho T''_1\rho_1( E^{(p)})
=(-1)^pv^{p(p-1)}E^{(p)}K^p$.
\end{proof}
A similar argument shows that if that $\mathcal F$ and $\mathcal E$ are finite dimensional $\bu$-modules
with $T''_1:\mathcal F^{T'_{-1}} \to \mathcal F$, $T''_1:\mathcal E^{T'_{-1}}\to
\mathcal E$, $\bu$-module homomorphisms and 
$\phi$ a $\rho$-invariant form on $\mathcal F\times \mathcal E$, then
\begin{equation}\label{L2}
\phi\circ(T''_1\otimes T''_1)\circ L^{-1}=\phi\circ
(T''_1\otimes T''_1) \circ (\mathcal  L^{-1}\otimes 1).
\end{equation}   
Observe that $\rho(\mathcal  L^{-1})=\mathcal  L ^{-1}$.

On the other hand
\begin{equation}\label{L3}
\phi\circ L=\phi\circ
(\mathcal  L\otimes 1).
\end{equation}   
Indeed
\begin{align*}
\phi\circ L(m\otimes n)&=\sum_{p\geq 0}(-1)^pv^{\frac{-p(p-1)}{2}}\{p\}
               \phi(F^{(p)}m, E^{(p)}*n) \\
&=\sum_{p\geq 0}(-1)^pv^{\frac{-p(p-1)}{2}}\{p\}\phi(\rho(\rho_1 E^{(p)})(F^{(p)}m),n)\\
&=\sum_{p\geq 0}v^{\frac{-3p(p-1)}{2}}\{p\}\phi( E^{(p)}K^{-p}F^{(p)}m),n).
\end{align*}
Similarly we have
$\rho(\mathcal L)=\mathcal L$.  Observe that $\mathcal L=\overline{\mathcal L^{-1}}$ where $\bar{\enspace}:\bu \to \bu$ is the automorphism of $\mathbb
Q$-algebras defined by
$\overline{E}=E$, $\overline{F}=F$,
$\overline{K}=K^{-1}$ and $\overline{fx}=\overline{f}\overline{x}$ for $f\in k(v)$, $x\in\bu$.  Here $\overline{v}=v^{-1}$ with
$\bar{\enspace}$ fixing $\mathbb Q$.  See \cite[3.1.12]{MR94m:17016} for more details. 
%

Moreover we can prove that $\mathcal  L^{-1}$ is an
intertwining map i.e. $\mathcal  L^{-1}:M^{T_1'}_\pi\to
M^{T_{-1}'}_\pi$ is a module  homomorphism for $M$ a highest weight module:

\begin{lem}\label{TandL} $T'_{-1}(u)\mathcal  L^{-1}=\mathcal  L^{-1}T_1'(u)$ as
operators on $M_\pi$ for all $u\in{_R\mathbf U}$.  From this we also get
$T'_1( u)\mathcal  L=\mathcal  LT_{-1}'( u)$ for all $u$.
\end{lem}
\begin{proof}  The proof follows from a direct calculation on generators of $\bu$.
\end{proof}

%

\section{Invariant Forms and Liftings} 

\subsection{} Elements of $\mathbb P(M,N)$ are called {\it invariant pairings} and for $M=N
$, set
$\mathbb P(M) =\mathbb P(M,M)$ and call the elements {\it invariant forms} on
$M$.

\begin{lem}  \begin{enumerate}
\item For $\ru$-modules $M,N$ and $\mathcal F$, 
\begin{equation}
\mathbb P(M,N,\mathcal F)\cong \hom_\ru (M\otimes_R N^{\rho_1},{_R\mathcal
F}^\brho)
\end{equation}
\item If $\mathcal F=R$ is the trivial $\ru$-module, then the condition $\phi\in\mathbb P(M,N)$, is equivalent to
\begin{equation}
\phi(u\,a,b)=\phi(a,\rho(u)b)
\end{equation}
for all $a\in M$, $b\in N$ and $u\in\bu$.

\end{enumerate}
\end{lem}

\begin{proof} The first statement follows from \cite[3.10.6]{MR96m:17029} and the second follows from a computation using the
generators of $\bu$.
\end{proof}

\begin{lem}
Let $M$ and $N$ be finite dimensional $\bu$-modules in
$\mathcal C'$ and   $\phi$ an invariant pairing. Then, for $m\in {_RM},n\in {_R N}$, $
\phi(\tiepp m,\tiempp n)=\phi(m,n)$.
\end{lem}

\begin{proof} We may assume $M$ and $N$ are irreducible with highest weight
$\nu$. Using the basis in (2.1), 
$h+j=\nu$, and the invariance of $\phi$ we get
\begin{align*}\phi(\tiepp F^{(j)}\eta,\tiempp F^{(j)}\eta)&=
\phi(F^{(h)}\eta,F^{(h)}\eta)
        =\phi(F^{(j)}\eta,F^{(j)}\eta)
\end{align*} where we have used \cite[Cor.
3.1.9]{MR94m:17016}, in the third line.
\end{proof}

Note that if $\phi(\eta,\eta)=1$, then for $0\leq j\leq \nu$, the calculation required for the
proof above shows that 
\begin{equation}\label{normalizedforms}
\phi(F^{(j)}\eta,F^{(j)}\eta)=v^{j^2-\nu j}{\qbinom \nu j}
\end{equation}

For the proof of some future results we must be explicit about the 
definition of the "$R$-matrix" $_f\mathcal R:A\otimes
\mathcal E\to \mathcal E\otimes A$:   Recall a
$\bu$-module
$M$ is said to be {\it integrable} if for any $m\in M$ and all $i\in I$, there exists a
positive integer
$N$ such that $E^{(n)}_im=0=F^{(n)}_im$ for all $n\geq N$ , and $M=\oplus
_{\l\in X}M^\l$ where for any
$\mu\in Y, \l\in X$ and $m\in M^\l$ one has
$K_\mu m=v^{\angles{\mu ,\l}}m$. Let $R$ denote a commutative algebra over the ring
$\mathbb Q(v)[T^{\pm 1}]$ (such as in the introduction) and let
$f:X\times X\to R$ be a function such that  
\begin{align*}
f(\zeta,\zeta'\pm i')&=f(\zeta,\zeta')v^{ \mp \angles{i,\zeta}(i\cdot i/2)}  \\
f(\zeta\pm i',\zeta')&=f(\zeta,\zeta')v^{ \mp \angles{i,\zeta'}(i\cdot i/2)}T^{ \mp(i\cdot i)/2}
\end{align*}
for all $\zeta,\zeta'\in X$ and all $i\in I$ (see [L, 32.1.3] or
[Ja, 3.15]). 

Such a function $f$ exists: Let $H$ denote a set of coset representatives of $X/\mathbb Z[I]$,  let $c:H\times H\to \mathbb Z$ denote an arbitrary function and set 
$$
f(h+\nu,h'+\nu'):=v^{c(h,h')-\sum_i\nu_i\langle i,h'\rangle (i\cdot i)/2-\sum_i\nu_i'\langle i,h\rangle (i\cdot i)/2-\nu\cdot \nu'}T^{-\sum \nu_i(i\cdot i)/2}
$$
for $h,h'\in H$ and $\nu,\nu'\in \mathbb Z[I]$.  

\begin{thm}{[L, 32.1.5], or [Ja, 3.14]} If $\mathcal E$ is an
integrable
$\ru $ module and $A\in \mathcal C_{R}$, then for each $f$ satisfying (3.1.1), there
exists an isomorphism $_f\mathcal R:A\otimes
\mathcal E\to \mathcal E\otimes A$. \end{thm}

The map $\tau:A\otimes B\to B\otimes A$ for any two modules $A$ and $B$ denotes the
twist map
$\tau(a\otimes b)=b\otimes a$.   Define $\prod _f:\in\End_R({_R\mathcal E}\otimes{_R\mathcal F}\otimes_R M)$ by 
$\prod_f(e\otimes e'\otimes  m)= f(\lambda,\lambda')e\otimes e'
\otimes  m$ for $m\in M^{\lambda'}$ and
$e\otimes e'\in (\mathcal E\otimes{_R\mathcal F})^{\lambda}$.  Lastly we define
$\chi\in
\End_R({_R\mathcal E\otimes{_R\mathcal F}}\otimes_R M)$ by
$$
\chi(e\otimes e'\otimes m)=
\sum_\nu\sum_{b,b'\in
\mathbf B_\nu}p_{b,b'}b^-(e\otimes e')\otimes {b'}^+m
$$ 
where $p_{b,b'}=p_{b',b}\in R$, and $\mathbf B_\nu$ is a subset of
$\mathfrak f$.   Then $_f\mathcal R$ is
defined to be equal to $\chi\circ \prod_f\circ \tau$.  The proof that it is an
$\bu$-module homomorphism is almost exactly the same as in [L, 32.1.5] or [Ja,
3.14], which the exception that one must take into account that $M$ is in the
category
$\mathcal C_R$ instead of $_R\mathcal C'$. 
(The above formula for 6.2.1 corrects
the corresponding 6.2.1 in the earlier version of [CE]. )

 In the case that $I=\{i\}$ and $X=Y=\mathbf Z$ with $i=1\in Y$, $i'=2\in X$, the linear map $L$ (see \eqnref{L}) coincides with $\chi$ (Lusztig uses the notation $\Theta$ where we use $\chi$ - see \cite[4.1.4]{MR94m:17016}).  Moreover since $X=\mathbb Z\supset 2\mathbb Z=\mathbf Z[i]$ there are only two cosets in $X/\mathbb Z[i]$.   In this case the function $f$ above can be defined by
\begin{equation}\label{fexample}
 f(h+\nu,h'+\nu'):=v^{-\sum_i\nu_i\langle i,h'\rangle (i\cdot i)/2-\sum_i\nu_i'\langle i,h\rangle (i\cdot i)/2-\nu\cdot \nu'}T^{-\sum \nu_i(i\cdot i)/2}
\end{equation}
for all $h,h'=0,1$ and $\nu,\nu'\in 2\mathbb Z$.  So for example $_f\mathcal R^{-1}$ is given by 
\begin{align}\label{commutativity} 
_f\mathcal R^{-1}(e\otimes m)&= \mathbf s\circ \Pi_f^{-1}\circ L^{-1}(e\otimes m)  \\
&=\sum_{n\geq 0} v^{n(n-1)/2}\{n\} \mathbf s\circ \Pi_f^{-1}F^{(n)}e\otimes E^{(n)}m \notag \\
&=\sum_{n\geq 0} f(\lambda-2n,\lambda'+2n)^{-1}v^{n(n-1)/2}\{n\}E^{(n)}m\otimes F^{(n)}e\notag \\
&=f(\lambda,\lambda')^{-1}\sum_{n\geq 0} v^{n(\lambda-\lambda')-2n^2+n(n-1)/2}T^{-n}\{n\}E^{(n)}m\otimes  F^{(n)}e\notag
\end{align}
for $e\in \mathcal E^{\lambda}$ and $m\in M^{\lambda'}$.

Let $\mathcal E$ and $\mathcal F$ be finite dimensional $\bu$-modules and 
$\tau: _R\mathcal E\otimes _R\mathcal F^{\rho_1}  \rightarrow \bu$, a
$\bu$-module homomorphism into $\bu$, where $\bu$ is a module under
the adjoint action.
Suppose $\phi$
is a pairing of $M$ and $N$.  Define $\psi_{\tau,\phi}$ to be  the
invariant pairing of $ M\otimes _R\mathcal E$ and $N\otimes _R\mathcal F$ defined
by
the formula, for $e\in \mathcal E ,f\in  \mathcal F,m\in M,$ and $n\in N$, 
\begin{equation}\label{inducedpairing}
\psi_{\tau,\phi} (m\otimes e,n \otimes f)=
\phi( m,\tau(e\otimes f)*n)\ ,  
\end{equation}
Here  ${\mathcal F}^{\rho_1}$ is a twist of the representation
$\mathcal F$ by $\rho_1$.  We call the pairing $\psi_{\tau,\phi}$ the
{\it pairing induced} by
$\tau$ and $\phi$. In the cases when $M,N$ and $\phi$ are fixed we write
$\psi_\tau$ in place of $\psi_{\tau,\phi}$ and say this pairing is induced
by $\tau$.

\subsection{} 
A result from \cite{MR96c:17022} shows that in the setting of Verma
modules the collection of maps $\tau$ is a natural set of parameters for
invariant forms. 

\begin{prop}
Suppose $\bu$ is of finite type and $\mathcal E$ and
$\mathcal F$ are finite dimensional $\bu$-modules. Let  $M$ be an
$_R\bu$-Verma module  and $\phi$ the Shapovalov form on $M$.
Then
every invariant pairing of
$M\otimes _R\mathcal E$ and $M\otimes _R\mathcal F$ is induced by $\phi$.
\end{prop}

\subsection{}

\begin{thm}[Lifting Theorem] 
Let $A$ and
$B$ be modules in $\mathcal C_{R}$ and $\phi\in {\mathbb P}_{\brho }(A,B)$. Then
$\phi$ uniquely determines an invariant form
$\phi_F\in  {\mathbb P}_{\brho}(A_F,B_F)$ which is determined by the following
properties:
\begin{enumerate}
\item $\phi_F$ vanishes on the subspaces $\i A\times B_F$ and
$A_F\times \i B$ .
\item For each $\mu \in \mathbb Z$ with $\mu+1= r\in \mathbb N$, and any vectors $a\in A$ and $b\in B$ both of weight $m+\e\mu$ with
$\e\in
\{1,s\}$ and  $E\ a=E\ b=0$ , 
\begin{equation}\label{lift}
\ \phi_F(F^{-1}a ,F^{-1} b)=  v^{-r+1}\frac{\i \e[T;0]}{  \i
\e[T;r-1]\ }\ \
\phi(a,b).
\end{equation}
\end{enumerate}
\end{thm}

\begin{prop}\label{invariance}  The form
$\phi_F$ induces an $\brho$-invariant bilinear map on
$A_\pi \times B_\pi$ which we denote by
$\phi_\pi$.
\end{prop}        

\subsection{} At times the subscript notation for lifted forms will be
inconvenient and so we shall also use the symbol $loc$ for the localization
of both forms and modules. We write $loc(\phi)$ and $loc(A)$ in place of
$\phi_F$ and $A_F$. 

For invariant forms we find that induction and localization commute in the
following sense.
\begin{prop}[{\cite[Prop 7.5]{MR96c:17022}}]  Suppose $A$ and $B$ are objects in
$\mathcal C_{R}$. Then $\phi_F$ and $\phi_\pi$ are $\ru$-invariant; i.e.        $\phi_F\  \in {\mathbb
P}_\brho(A_F,B_F)\  ,\
\phi_\pi\in {\mathbb P}_\brho(A_\pi,B_\pi)\ .$
\end{prop}

\section{Quantum Clebsch-Gordan decomposition}\label{QCG decomposition}

\subsection{Bases and Symmetries} For
$m\in\mathbb Z$, let
$\mathcal F_m$\label{fm} denote the finite dimensional irreducible module of
highest weight
$v^m$ with highest weight vector $u^{(m)}$.  For $k$ any non-negative
integer set
$u^{(m)}_k=F^{(k)}u^{(m)}$ and
$u^{(m)}_{-1}=0$.

In particular 
$$
T_1''(u^{(m)}_j)=(-1)^{m-j}v^{(m-j)(j+1)}u^{(m)}_{m-j}
$$
and
\begin{equation}\label{irred}
K^pu^{(m)}_j=v^{p(m-2j)}u^{(m)}_j,\enspace  F^{(p)}u^{(m)}_j=
{\qbinom {p+j} j}u^{(m)}_{j+p},\enspace
 E^{(p)}u^{(m)}_j=\qbinom{m+p-j}{
p} u^{(m)}_{j-p}.
\end{equation}

\begin{lem}[Clebsch-Gordan, \cite{MR96e:17041},
\cite{MR90m:17022}] \label{CG} For any two non-negative integers
$m$ and $n$, there is an isomorphism of
$\bu$-modules
$$
\mathcal F_{m+n}\oplus \mathcal
F_{m+n-2}\oplus \cdots \oplus \mathcal F_{|m-n|}\cong\mathcal F_m\otimes
\mathcal F_n.
$$
Moreover the isomorphism may be defined on highest weight vectors by
\begin{equation}
\varPhi(u^{(m+n-2p)})=\sum_{k=0}^{p}(-1)^kv^{(k-p)(m-p-k+1)}\frac{[n-p+k]_{(k)}}{[m]_{(k)}}u_k^{(m)}\otimes u_{p-k}^{(n)}
\end{equation}
where $0\leq p\leq \min\{m,n\}$.
\end{lem}

\begin{proof}
First note $|m-n|\leq m+n-2p\leq m+n$ so that $0\leq p\leq\min\{m,n\}$.
 The element $\varPhi(u^{(m+n-2p)})$ certainly has the right weight and it is straightforward to check
 that $E\varPhi(u^{(m+n-2p)})=0$.
\end{proof}

\begin{cor} For any two non-negative integers $m$, $n$, and $p\leq\min\{m,n\}$,
\begin{align*}
\varPhi(u^{(m+n-2p)}_{m+n-2p})
&=v^{p(n- m)}\sum_{k=0}^{p}(-1)^{p+k}v^{-k(n+k-2p+1)}\frac{[n+k-p]_{(k)}}{[m]_{(k)}}
  u_{m-k}^{(m)}\otimes u_{n+k-p}^{(n)}.
\end{align*}

\end{cor}

\begin{lem}  The map $\varphi:\mathcal F_m^{\rho_1}\to \mathcal F_m$ given
by $\varphi(u^{(m)}_k)=(-v)^{-k}u^{(m)}_{m-k}$ is an isomorphism.
\end{lem}
\begin{proof} This is a straightforward calculation on weight vectors and generators of $\bu$. 

\end{proof}

\begin{cor}\label{CG2} Let
$m$ and $n$ be two non-negative integers. Then there is an isomorphism of
$\bu$-modules
$$
\mathcal F_{m+n}\oplus \mathcal
F_{m+n-2}\oplus \cdots \oplus \mathcal F_{|m-n|}\cong\mathcal F_m\otimes
\mathcal F_n^{\rho_1}.
$$
Moreover for $m\geq n$, this isomorphism can defined on highest weight vectors by
\begin{align*}
\Phi(u^{(m+n-2p)})&=\sum_{k=0}^{p}(-1)^{n-p}\frac{[n-p+k]![m-k]!}{[n-p]![m]!}
v^{(k-p)( 2 + m)+p^2 -k^2  + n}u_k^{(m)}\otimes
u_{n-p+k}^{(n)}
\end{align*}

(the action on the second factor $u^{(n)}_l$ is twisted by the
automorphism $\rho_1$).
\end{cor}

\begin{proof}  This follows from the isomorphism $\phi:\mathcal
F_n^{\rho_1}\to\mathcal F_n$ which sends
$u^{(n)}_{n+k-p}$ to $(-v)^{p-n-k}u^{(n)}_{p-k}$.
\end{proof}

Since we have two basis $\{u_i^{(m)}\otimes u_j^{(n)}|\,0\leq i\leq
m,\enspace 0\leq j\leq n\}$, and $\{u_k^{(m+n-2p)}|\,0\leq p\leq
n,\enspace 0\leq k\leq m+n-2p\}$ of $\mathcal F_m\otimes \mathcal
F_n$, we can relate them by the {\it quantum Clebsch-Gordan
coefficients} or {\it quantum $3j$-symbols};
$$
u_k^{(m+n-2p)}=
\sum_{0\leq i\leq m,0\leq j\leq n,i+j=p+k}
\left[\begin{matrix} m & n & m+n-2p
\\ i & j & k \end{matrix}\right]u_i^{(m)}\otimes u_j^{(n)}
$$

Consider now the $\rho$-invariant forms \eqnref{normalizedforms} on
$\mathcal F_m$ and
$\mathcal F_n$, both denoted by $(,)$, normalized so that their
highest weight vectors have norm $1$.  Define the symmetric invariant
bilinear form on $\langle,\rangle$ on $\mathcal F_m\otimes \mathcal
F_n$ given by the tensor product of the two forms (the resulting
pairing is $\rho$-invariant). Assume that the forms on $\mathcal F_m$ and $\mathcal F_n$ are normalized so that their highest weight vectors
$u^{(m)}$ and $u^{(n)}$ have norm $1$. In this case
\begin{align*}
\langle u^{(m+n-2p)},u^{(m+n-2p)}\rangle
=
v^{p\,\left(2\,p - 2\,m
-1\right)}\frac{[n]![m+n-p+1]![m-p]!}{[m]![p]![m+n-2p+1]![n-p]!}.
\end{align*}
where we have used formula \eqnref{Ma2} for $p\leq \min\{m,n\}$.

The same proof that gave us \eqnref{normalizedforms} now implies for $p\leq\min\{m,n\}$,
\begin{align*}
||u^{(m+n-2p)}_k||^2:=v^{p\,\left(2\,p - 2\,m -1\right)-(m+n-2p-k) k}
\frac{\left[\begin{matrix} n \\  p
\end{matrix}\right]\left[\begin{matrix} m+n-p+1
 \\  p \end{matrix}\right]\left[\begin{matrix} m+n-2p \\ k\end{matrix}\right]}{\left[\begin{matrix} m \\  p
\end{matrix}\right]}
\end{align*}

\begin{prop}[\cite{MR90m:17022}]
\begin{enumerate}[(i).]
\item The basis $\{u_k^{(m+n-2p)}\}$
of \hbox{$\mathcal F_m\otimes{\mathcal F_n}$} is orthogonal.
\item For $0\leq i\leq m$, and $0\leq j\leq n$, 
\begin{align*}
&u_i^{(m)}\otimes u_j^{(n)} \\
&=v^{mj+ni-2ij} 
   {\qbinom m i}{\qbinom n j} \\
&\quad \times \sum_{p=0}^{\min\{m,n\}}\frac{[i+j-p]![m+n-i-j-p]!}{v^{p(m+n-p)}||u^{(m+n-2p)}||^{2}[m+n-2p]!}
  \left[\begin{matrix} m & n & m+n-2p
  \\ i & j & i+j-p \end{matrix}\right]u_{i+j-p}^{(m+n-2p)}.
\end{align*}
\end{enumerate}
\end{prop}

In \hbox{$\mathcal F_m\otimes{\mathcal F_n}^{\rho_1}$} 
(recall $\phi(u^{(n)}_j)=(-v)^{-j} u^{(n)}_{n-j}$ )
\begin{align}
&u^{(m)}_i\otimes u_j^{(n)} \\
&=(-1)^jv^{m(n-j)+ni-2i(n-j)-j} 
   {\qbinom m i}{\qbinom{n}{n-j}}\notag \\
&\quad \times \sum_{p=0}^{\min\{m,n\}}\frac{v^{-p(m+n-p)}}{||u^{(m+n-2p)}||^{2}} \notag\\
&\qquad\times\frac{[n+i-j-p]![m-i+j-p]!}{[m+n-2p]!}
  \left[\begin{matrix} m & n & m+n-2p
  \\ i & {n-j} & n+i-j-p \end{matrix}\right]u_{n+i-j-p}^{(m+n-2p)}.\notag
\end{align}

\begin{proof}  First we observe that due to fact that different weight
spaces having different weights we get
$$
\langle u_k^{(m+n-2p)},u_l^{(m+n-2p)}\rangle=0,
$$
for $k\neq l$.  Now suppose $p\neq q$.  Then the pairing $(,)$ on
$\mathcal F_m\otimes \mathcal F_n$ induces a module homomorphism from
$\mathcal F^{(m+n-2p)}\to ((\mathcal F^{(m+n-2q)})^\rho)^*$.  Since
these two irreducible modules are not isomorphic for $p\neq q$, we must
have
$$
\langle u_k^{(m+n-2p)},u_l^{(m+n-2q)}\rangle=0,
$$
for all $k$ and $l$.

The second formula follows from the fact that if we write 
$$
u_i^{(m)}\otimes u_j^{(n)}=\sum_{p=\frac{\left|m-n\right|}{2}}^{\frac{m+n}{2}}c^{ij}_p\,u_{i+j-p}^{(m+n-2p)},
$$
then $k=i+j-p$,
\begin{align*}
c^{ij}_p\langle u_k^{(m+n-2p)},u_k^{(m+n-2p)}\rangle 
&=\langle u_i^{(m)}\otimes u_j^{(m)},u_k^{(m+n-2p)}\rangle \\
&=\left[\begin{matrix} m & n & m+n-2p \\ i & j & k \end{matrix}\right](u_i^{(m)},u_i^{(m)})(u_j^{(n)}, u_j^{(n)})\\
&=v^{i(i-m)+j(j-n)}\left[\begin{matrix} m & n & m+n-2p \\ i & j & k \end{matrix}\right]{\qbinom m i}{\qbinom n j}
\end{align*}
\end{proof}

\begin{lem} For $m$ and $p$ non-negative integers, one has
\begin{equation}  
\mathcal L^{-1}u^{(m)}_{m-p}=v^{2p(p-1-m)}u^{(m)}_{m-p}
\end{equation}
and
\begin{equation}  
\mathcal Lu^{(m)}_{m-p}=v^{-2p(p-1-m)}u^{(m)}_{m-p}.
\end{equation}
These two equations explain the labeling of $\mathcal L^{-1}$ and $\mathcal L$ i.e. they really are inverses when restricted to finite dimensional
$X$-admissible modules. 
\end{lem}
\begin{proof}  This follows directly from \lemref{TandL} 
\end{proof}
\begin{cor} Suppose $m$, $n$, $k$ and $p$ are non-negative integers with $c^{m,n}_{i,j}\in k(v)$ and
$$ 
\bar u=\sum_{i+j=p+k}c^{m,n}_{i,j}u^{(m)}_i\otimes u^{(n)}_j.
$$
Then
\begin{equation}  
T''_1\mathcal L^{-1}\bar u=\sum_{i+j=p+k}d^{m,n}_{i,j}u^{(m)}_i\otimes u^{(n)}_j,
\end{equation}
where 
\begin{align}\label{key}
d^{m,n}_{r,s}&=(-1)^{r+s}v^{r(r - m-1 )+s(s - n-1)} \\ 
&\qquad \times\sum_{p\geq 0}(-1)^pv^{-\frac{p(p-1)}{2}+ p(2p+ m - n - 2r + 2s )}\left\{p\right\}
c^{m,n}_{m-r+p,n-s-p}\qbinom{r}{p}\qbinom{n-s}{p}\notag
\end{align}

\end{cor}
\begin{proof}  Due to the invariance of the form $\langle\,,\,\rangle$ we get
\begin{align*}
&d_{r,s}^{m,n}v^{r(r-m)+s(s-n)}\qbinom{m}{r}\qbinom{n}{s}=d_{r,s}^{m,n}(u^{(m)}_r,u^{(m)}_r)(u^{(n)}_s,u^{(n)}_s) \\
&=\langle T''_1\mathcal L^{-1}\bar u, u^{(m)}_r\otimes u^{(n)}_s\rangle  
=\langle T''_1\mathcal L^{-1}\bar u, T''_1T'_{-1}(u^{(m)}_r\otimes u^{(n)}_s)\rangle   \\
&=\langle L(\mathcal L^{-1}\bar u\otimes T'_{-1}(u^{(m)}_r\otimes u^{(n)}_s)\rangle  \qquad \text{by \eqnref{L2}} \\
&=\langle \bar u,T'_{-1}(u^{(m)}_r\otimes u^{(n)}_s)\rangle \qquad \text{by \eqnref{L3}}   \\
&=\langle \bar u,L(T'_{-1}(u^{(m)}_r)\otimes T'_{-1}(u^{(n)}_s))\rangle   \\
&=(-1)^{r+s}v^{r(2r-2m-1)+s(2s-2n-1)} \\ 
&\qquad \times\sum_{p\geq 0}(-1)^pv^{-\frac{p(p-1)}{2}+ p(2p+ m - n - 2r + 2s )}\left\{p\right\}
c^{m,n}_{m-r+p,n-s-p}\qbinom{r}{p}\qbinom{n-s}{p}\qbinom{m}{r}\qbinom{n}{s}.
\end{align*}

\end{proof}
\section{Basis and the Intertwining map $\mathcal L$}
\subsection{A Basis} For $s\geq 1$,
and any lowest weight vector $\eta$ of weight $Tv^{\lambda+\rho}$, set
\begin{equation}\label{fminusketa}
F^{(-k)}\eta:=T'_{-1}(F^{(k)})\eta
=v^{k(k-1)}F^{-k}K^{k}[K;-1]^{(k)}\eta
=v^{k(\lambda+k)}T^k[T;\lambda]^{(k)}F^{-k}\eta
\end{equation}

as
$$
T'_{-1}(E^{(k)})=(-1)^kv^{-k(k-1)}K^{-k}F^{(k)},\quad
\quad T'_{-1}(F^{(k)})=(-1)^kv^{k(k-1)}E^{(k)}K^{k}.
$$
\begin{lem}\label{binomiallemma} Suppose $r,s\in\mathbb Z$, $s>
r\geq 0$, $\zeta$ is a highest weight vector of weight
$Tv^{\lambda-\rho}$ and $\eta$ is a lowest weight vector of
weight $Tv^{\lambda+\rho}$.  Then
\begin{equation}\label{oldandnew}
E^{(r)}F^{(s)}\zeta
=\qbinom{T;\lambda-\rho+r-s}{ r}F^{(s-r)}\zeta,
\quad F^{(r)}F^{(s)}\zeta
=\qbinom{r+s}{ r}F^{(s+r)}\zeta, 
\end{equation}
and
\begin{align}
F^{(r)}F^{(-s)}\eta
&=v^{r(\lambda+2s-r)}T^r
 \qbinom{T;\lambda+s}{ r} F^{(r-s)}\eta, \\
E^{(r)}F^{(-s)}\eta
&=(-1)^rv^{-r(\lambda+r+2s)}T^{-r}
 \qbinom{r+s}{ r}
  F^{(-r-s)}\eta, 
\end{align}

\end{lem}

$n=\lambda-\rho$.

Define indexing sets $ I_\l$ and $ I_{-\l}$ by $ I_\l =
\lbrace n-2,n-4,...\rbrace$, $ I_{-\l} = \lbrace -n,-n-2,...\rbrace $ where $n=\lambda+1$.  
One should compare the previous result with
\begin{lem}[\cite{MR96c:17022}, 2.2]\label{locstructure} Now for
integers $j \in  I_\l$ (resp.
$ I_{-\l}$) set $ k_j = {\frac{n-2-j}{ 2}}$ and $l_j = {\frac{-n-j}{2}}
$ and define basis vectors for $M(m+\l)$ and $M(m-\l)$ by 
$\vj = F ^{k_j} \otimes 1_{m+\ll}$ and $\vsj = F ^{l_j}\otimes 1_{m-\l
-\rho}$. The action of $\ru$ is given by 
$$ 
K\vj = Tv^{\l-1-2k_j}\vj
\quad ,\quad F \vj = w_{\l,j-2}\ \  , 
$$ 
$$ K\vsj =Tv^{-\l-1-2l_j}\vsj \quad ,\quad  F
\vsj = w_{-\l,j-2} \ , 
$$ 
$$ E\vj =  [k_j][T; -l_j] w_{\l,j+2}\quad ,\quad    E\vsj = [l_j][T; -k_j]
w_{-\l,j+2}. 
$$

\end{lem}

\subsection{} The
articles \cite{MR96c:17037} and
\cite{MR96c:17022} study noncommutative localization of highest
weight modules. This article may be viewed as an extension of what
was begun there. For any $\bu$-module $A$ let $A_F$ denote the
localization of
$A$ with respect to the multiplicative set in $\bu$ generated by $F$. If
$F$ acts without torsion on $A$ (we shall assume this throughout) then $A$
injects into $A_F$ and we have the short exact sequence of $\bu$-modules:
$0\ra A\ra A_F\ra A_\pi \ra 0$. 
\begin{lem}[Mackey] \label{Mackey}  Let $M=\ru\otimes _{_RB}R_{m+\lambda-\rho}$ be a Verma module of highest weight
$m+\lambda-\rho$ (see \cite[1.7]{MR96c:17022}) and suppose that $\mathcal E$ is a $X$-admissible (hence a free $R$-module of
finite rank). Then $M_\pi\otimes_R\mathcal E$ is generated as an $\ru$-module by $Rm\otimes \mathcal E$ where $m$ is a
lowest weight vector in $M_\pi$.  
\end{lem}
\begin{proof}  First note that $M\otimes_R \mathcal E$ can be identified with its image in $M_F\otimes R\mathcal E$ due to the
fact that it is torsion free with respect to the action of $F$.  The Lemma then follows from the series of isomorphisms:
\begin{align*}
M_\pi\otimes_R\mathcal E&\cong (M_F\otimes_R\mathcal E)/(M\otimes_R\mathcal E)\cong (M\otimes_R\mathcal E)_F/(M\otimes_R\mathcal E) \\
			&\cong (\ru\otimes_{_RB}(R_{m+\lambda-\rho}\otimes_R\mathcal E))_F/
   (\ru\otimes_{_RB}(R_{m+\lambda-\rho}\otimes_R\mathcal E)) \\
			&= (\ru\otimes_{_RB}(R_{m+\lambda-\rho}\otimes_R\mathcal E))_\pi.
\end{align*}
The first isomorphism follows as tensoring the short exact sequence 
$$
0\to M\to M_F\to M_\pi\to 0
$$
with $\mathcal E$ (a finite rank free $R$-module), leads to a short exact sequence.  The second isomorphism comes from
\cite[Theorem 3.9]{MR96c:17022}) where it is implemented by $1\otimes e\mapsto 1\otimes 1\otimes e$ where $e\in\mathcal E$. 
The third isomorphism comes from Mackey's Isomorphism Theorem and it is implemented by $1\otimes 1\otimes e\mapsto
1\otimes 1\otimes e$.
\end{proof}

\subsection{Highest and Lowest Weight Decompositions.}

\begin{lem}\label{highestlowestwtvectors}  Suppose $n-2p=r$ with $n$ and $r$ nonnegative integers.
 Then 
 \begin{equation}
w_{r,r-1}=\sum_{0\le k\le p}\frac{[n-p+k]_{(k)}v^{-k^2}}{T^{-k}
[T^{-1};k]_{(k)}} F^{(k)}w_{0,-1}\otimes u^{(n)}_{p-k}.
\end{equation}
is a highest weight vector of weight $Tv^{r-1}$ in $M\otimes \mathcal F$
and 
\begin{equation}
w_{-r,-r-1}=\sum_{0\le k\le n-p}\frac{[p+k]_{(k)}v^{-k^2}}{T^{-k}
[T^{-1};k]_{(k)}} F^{(k)}w_{0,-1}\otimes u^{(n)}_{n-p-k}.
\end{equation}
is a highest weight vector of weight $Tv^{-r-1}$ in $M\otimes \mathcal F$.  In $M_\pi\otimes \mathcal F$,
 \begin{equation}
m_{r,r+1}=\sum_{0\le k\le n-p}(-1)^k\frac{[p+k]_{(k)}v^{-k(-r+2k+1)}}{T^k
[T;k]_{(k)}} F^{(-k)}m_{0,1}\otimes u^{(n)}_{p+k},
\end{equation}
is a lowest weight vector of lowest weight $Tv^{r+1}$.  Moreover
 \begin{equation}
m_{-r,-r+1}=\sum_{0\le k\le p}(-1)^k\frac{[n-p+k]_{(k)}v^{-k(r+2k+1)}}{T^k[T;k]_{(k)}}F^{(-k)}m_{0,1}\otimes u^{(n)}_{n-p+k},
\end{equation}
is a lowest weight vector in $M_\pi\otimes \mathcal F$ of
lowest weight $Tv^{-r+1}$
\end{lem}

Observe that the second identity can be obtained from the first by replacing $p$ with $n-p$.
This implies that if we set
 \begin{equation*}\label{sigma}
\sigma=\sigma(\epsilon,n,p)=\begin{cases} p &\quad\text{ if }\quad \epsilon =1, \\
  n-p &\quad \text{if} \quad \epsilon =-1,\end{cases}
\end{equation*}
and
 \begin{align*}
a_{\epsilon, n, p}(k)=(-1)^k\frac{[\sigma+k]_{(k)}}{T^k[T;k]_{(k)}},
\end{align*}
then
 \begin{equation}
m_{\epsilon r,\epsilon r+1}=\sum_{0\le k\le n-\sigma}v^{-k(-\epsilon r+2k+1)}a_{\epsilon, n, p}(k)F^{(-k)}m_{0,1}\otimes
u^{(n)}_{\sigma+k},
\end{equation}

\begin{proof}  The first equality follows using \eqnref{irred} and \lemref{binomiallemma}.

On the other hand a direct calculation shows $Fm_{r,r+1}=0$.
Using 
\begin{align*}
m_{-r,-r+1}&=\sum_{0\le k\le p}A_k F^{(-k)}m_{0,1}\otimes u^{(n)}_{n+k-p}
\end{align*}
and $Fm_{-r,-r+1}=0$ we get
$$
A_k =(-1)^k\frac{v^{-k(r+2k+1)}[n+k-p]_{(k)}}{T^k[T;k]_{(k)}}A_0.
$$
\end{proof}
Note that by the definition of $\mathcal L$, one has
\begin{align}
\mathcal L(m_{r,r+1})=m_{r,r+1},\quad 
\mathcal L(m_{-r,-r+1})=m_{-r,-r+1}.
\end{align}
\begin{cor}
Suppose $n-2p=r$ with $n$ and $r$ nonnegative integers. If we let $m_{r,r+1}$ (resp. $m_{-r,-r+1}$) denote a lowest weight vector in $M_\pi\otimes
\mathcal F$ of lowest weight $Tv^{r+1}$ (resp.   $Tv^{-r+1}$),  then
 \begin{align}
L^{-1}(m_{r,r+1})&= v^{-2(n-p)(p+1)}\sum_{s=0}^{n-p}(-1)^s\frac{v^{s(1-r)}[p+s]_{(s)}}{T^s[T,s]_{(s)}}F^{(-s)}m_{0,1}\otimes
  u^{(n)}_{p+s},\label{needlater}  \\  
L^{-1}(m_{-r,-r+1})&=v^{-2p(n-p+1)}\sum_{s=0}^{p}\frac{v^{s(1+r)}[n-p+s]_{(s)}}{T^s[T,s]_{(s)}}F^{(-s)}m_{0,1}\otimes
u^{(n)}_{n-p+s}.
\end{align}
\end{cor}
As in the previous lemma the second identity can be obtained from the first by replacing $p$ with $n-p$. 
Using the above definition of $\sigma$ and $a_{\epsilon, n, p}(k)$ we get
 \begin{align*}
L^{-1}(m_{\epsilon r,\epsilon  r+1})&= v^{-2(n-\sigma)(\sigma+1)}\sum_{s=0}^{n-\sigma}(-1)^sv^{s(1-\epsilon
r)}a_{\epsilon, n, p}(k)F^{(-s)}m_{0,1}\otimes u^{(n)}_{\sigma+s}.
\end{align*}
\begin{proof}  Define coefficients $B_s$ and $C_s$ through
 \begin{align*}
L^{-1}(m_{r,r+1})  
&= \sum_{s=0}^{n-p}B_sF^{(-s)}m_{0,1}\otimes u^{(n)}_{p+s},\quad \text{and}  \\  
L^{-1}(m_{-r,-r+1})  
&= \sum_{s=0}^{p}C_sF^{(-s)}m_{0,1}\otimes u^{(n)}_{s+n-p}.
\end{align*} 
Now 
\begin{align*}
0&=L^{-1}(Fm_{\epsilon r,\epsilon r+1})
=-(K\otimes K)(T'_{-1}(E)\otimes 1+K^{-1}\otimes T'_{-1}(E))L^{-1}(m_{\epsilon r,\epsilon r+1})
\end{align*}
Thus
 \begin{align*}
0&= -\sum_{s=1}^{n-p}\left(B_s[T;s]+B_{s-1}T^{-1}v^{1-n+2p}[p+s]\right)F^{(1-s)}m_{0,1}\otimes u^{(n)}_{p+s}    ,
\end{align*} 
and
 \begin{align*}
0&= -\sum_{s=1}^{p}\left(C_s[T;s]+C_{s-1}T^{-1}v^{n-2p+1}[n-p+s]\right)F^{(1-s)}m_{0,1}\otimes u^{(n)}_{s+n-p}.
\end{align*}

Hence
\begin{align*}
B_s&=-B_{s-1}\frac{v^{1-n+2p}[p+s]}{T[T,s]}=(-1)^s\frac{v^{s(1-r)}[p+s]_{(s)}}{T^s[T,s]_{(s)}}B_{0} \\
C_s&=-C_{s-1}\frac{v^{n-2p+1}[n-p+s]}{T[T,s]}=(-1)^s\frac{v^{s(r+1)}[n-p+s]_{(s)}}{T^s[T,s]_{(s)}}C_{0}.
\end{align*}
In particular $B_{n-p}=(-1)^{n-p}\frac{v^{(n-p)(1-r)}[n]_{(n-p)}}{T^{n-p}[T;n-p]_{(n-p)}}B_{0}$.  On the other hand
\begin{align*}
&L^{-1}(m_{r,r+1})=
  \sum_{0\le k\le n-p}(-1)^k\frac{[p+k]_{(k)}v^{k(r-2k-1)}}{T^k
  [T;k]_{(k)}} L^{-1}\left(F^{(-k)}m_{0,1}\otimes u^{(n)}_{p+k}\right) \\
&= \sum_{s=0}^{n-p}\sum_{k=s}^{n-p}(-1)^k\frac{[p+k]_{(k)}}{
[T;s]_{(s)}[k-s]!}\qbinom{n-p-s}{k-s}v^{k(r-2k-1)+\frac{(k-s)(k-s-1)}{2}+(k-s)(k+s)}\\
&\hskip 100pt \times \{k-s\}T^{-s}F^{(-s)}m_{0,1}\otimes u^{(n)}_{p+s}. 
\end{align*}
A very similar calculation shows
\begin{align*}
L^{-1}&(m_{-r,-r+1}) \\
&= \sum_{s=0}^{p}\sum_{k=s}^{p}(-1)^k\frac{[n-p+k]_{(k)}}{
[T;s]_{(s)}[k-s]!}\qbinom{p-s}{k-s}v^{-k(r+2k+1)+\frac{(k-s)(k-s-1)}{2}+(k-s)(k+s)}\\
&\hskip 100pt \times \{k-s\}T^{-s}F^{(-s)}m_{0,1}\otimes
u^{(n)}_{n-p+s}  \\  
\end{align*}

When $s=n-p$ we get $F^{(-s)}m_{0,1}\otimes
u^{(n)}_{p+s}$ has a coefficient
$$
B_{n-p}=(-1)^{n-p}\frac{[n]_{(n-p)}v^{(n-p)(-n-1)}}{T^{n-p}
[T;n-p]_{(n-p)}}=(-1)^{n-p}\frac{v^{(n-p)(1-r)}[n]_{(n-p)}}{T^{n-p}[T;n-p]_{(n-p)}}B_{0}
$$
Thus $B_0=v^{(n-p)(r-n-2)}=v^{-2(n-p)(p+1)}$, 
$$
B_s=(-1)^s\frac{v^{s(1-r)-2(n-p)(p+1)}[p+s]_{(s)}}{T^s[T,s]_{(s)}},
$$
Similarly setting $s=p$, we get
$$
C_p=(-1)^p\frac{[n]_{(p)}}{T^p[T;p]_{(p)}}v^{-p(r+2p+1)}=(-1)^p\frac{[n]_{(p)}}{T^p[T;p]_{(p)}}v^{p(r+1)}C_{0}
$$
so that $C_0=v^{-2p(r+p+1)}=v^{-2p(n-p+1)}$.

This proves the two identities.
\end{proof}
\begin{cor}\label{link}  For $0\leq s\leq n-p$,
\begin{align}
&v^{s(1-r)+2(p-n)(p+1)}[p+s]_{(s)} \\
&=\sum_{k=s}^{n-p}(-1)^{k+s}\frac{[p+k]_{(k)}}{[k-s]!}\qbinom{n-p-s}{k-s}v^{k(r-2k-1)+\frac{(k - s)(3k + s-1) }{2}}\{k-s\}.\notag
\end{align}
\end{cor}

\section{Maps into the Harmonics.}  
\subsection{Harmonics}  We know from \cite{MR95e:17017} that
$F(\bu)\cong \mathcal H\otimes Z(\bu)$ where $\mathcal
H=\oplus_{n\in\mathbb N}\mathcal L_{2n}$\label{mathcalL2n} is the space of {\it
harmonics}, and
$\mathcal L_{2n}\cong
\mathcal F_{2n}$.

For $m$, $n$, $r$ integers with $m+n$ even, $|m-n|\leq 2r\leq m+n$, we let
$\beta^{m,n}_{2r}:\mathcal F_m\otimes\mathcal F_n^{\rho_1}\to
\mathcal H$, be the
$\bu$-module homomorphism determined by 
\begin{equation}
\label{betamn}
\beta^{m,n}_{2r}(u^{(m+n-2q)})=\delta_{2r,m+n-2q}E^{(r)}
K^{-r}.
\end{equation}
so that $\im \beta^{m,n}_{2r}=\mathcal L_{2r}$.

\begin{prop}\label{betaonabasis}Suppose $\eta$ is a lowest weight
vector of weight
$Tv^{\lambda+\rho}$, 
and
$\max\{|\frac{m-n}{2}|,|\frac{m-n}{2}-i+j|\}\leq r\leq \frac{m+n}{2}$ 
with $0\leq i\leq m$ and $0\leq j\leq n$.    Then
\begin{align*}
\rho_1\Big(\beta^{m,n}_{2r}(&u^{(m)}_i\otimes u^{(n)}_{j})
 \Big) F^{(-c)}\eta\\
&=(-1)^jv^{m(n-j)+ni-2i(n-j)-j+r^2-\frac{(m+n)^2}{4}} 
   {\qbinom m i}{\qbinom{n}{j}} \\ 
&\quad \times 
  \frac{\left[\begin{matrix} m & n & 2r
  \\ i & {n-j} & \frac{n-m}{2}+r+i-j \end{matrix}\right]}{||u^{(2r)}||^2\left[\begin{matrix} 2r
  \\  \frac{n-m}{2}+r+i-j \end{matrix}\right]}
(v^{2\lambda+2+i-j+\frac{n-m}{2}+4c}T^2)^{j-i+\frac{m-n}{2}}\\ 
&\quad \times 
 \sum_{l=0}^{i-j+\frac{n-m}{2}+r} (-1)^{r-l}
  v^{l\left(r+j-i+\frac{m-n}{2}+1 \right) } \\
&\hskip 50pt \times \qbinom{i-j+\frac{n-m}{2}+c}{ i-j+\frac{n-m}{2}+r-l}
 \qbinom{T;\lambda+l+c}{ r}
 \qbinom{l+c}{ l}
  F^{(j-i+\frac{m-n}{2}-c)}\eta \end{align*}
and
\begin{align*}
\rho_1\Big(T_{-1}'&\beta^{m,n}_{2r}(u^{(m)}_i\otimes u^{(n)}_{j})
 \Big) F^{(-c)}\eta\\
&=(-1)^jv^{m(n-j)+ni-2i(n-j)-j+r^2-\frac{(m+n)^2}{4}} 
   {\qbinom m i}{\qbinom{n}{j}} \\ 
&\quad \times ||u^{(2r)}||^{-2}
  \frac{\left[\begin{matrix} m & n & 2r
  \\ i & {n-j} & \frac{n-m}{2}+r+i-j \end{matrix}\right]}{\left[\begin{matrix} 2r
  \\  \frac{n-m}{2}+r+i-j \end{matrix}\right]}  (v^{ 2\,c
+\lambda}T)^{i-j+\frac{n-m}{2}}\\ 
&\quad\times  \sum_{l=0}^{i-j+\frac{n-m}{2}+r}(-1)^{r-l}
  v^{  l (r +j-i+\frac{m-n}{2}+1) }\\
&\hskip 50pt \times 
    \left[\begin{matrix}T;\lambda+c \\ l \end{matrix}\right]
    \left[\begin{matrix}T;\lambda+r+c-l \\ i-j+\frac{n-m}{2}+r-l \end{matrix}\right]
    \left[\begin{matrix}r+c-l \\ r \end{matrix}\right]
  F^{(i-j+\frac{n-m}{2}-c)}\eta  
\end{align*}

\end{prop}

\begin{proof}
From \cite{MR96m:17029}, 4.18(5), we have for any $a,b\in\mathbb N$,
$$
\text{ad}\,F^{(a)}(E^{(b)}K^{-b})=\sum_{m=0}^a(-1)^{a-m}
  v^{- \left( a-1 \right) \left(  a-m \right)}
  F^{(m)} E^{(b)}K^{-b}F^{(a-m)}K^{a}.
$$
This implies 
\color{black}

\begin{align*}
\rho_1&\left(\text{ad}\,F^{(a)}(E^{(b)}K^{-b})\right) F^{(-c)}\eta \\
&=\sum_{m=0}^a(-1)^{a-m}
  v^{- \left( a-1 \right) \left(  a-m \right)}
 \rho_1\left( F^{(m)} E^{(b)}K^{-b}F^{(a-m)}K^{a}\right) F^{(-c)}\eta \\
&=(v^{2\lambda+2+a-b+4c}T^2)^{b-a}
 \sum_m (-1)^{b-m}
  v^{m\left(2b-a+1 \right) }\\
&\hskip 50pt \times 
 \qbinom{a-b+c}{ a-m}
 \qbinom{T;\lambda+m+c}{ b}
 \qbinom{m+c}{ m}
  F^{(b-c-a)}\eta ,
\end{align*}

and since $\rho_1\circ T_{-1}'=T_{-1}'\circ \rho_1$, we get
\begin{align*}
\rho_1&\left(T_{-1}'\text{ad}\,F^{(a)}(E^{(b)}K^{-b})\right)
F^{(-c)}\eta \\
&=\sum_{m=0}^a(-1)^{a-m}
  v^{- \left( a-1 \right) \left(  a-m \right)}
 T_{-1}'\rho_1\left( F^{(m)} E^{(b)}K^{-b}F^{(a-m)}K^{a}\right)
  F^{(-c)}\eta \\ 
&=
 (-v^{ 2\,c +\lambda}T)^{a-b}\sum_{m}(-1)^{a-m}
  v^{  m (2b -a+1) }  \\
&\hskip 50pt \times 
    \left[\begin{matrix}T;\lambda+c \\ m \end{matrix}\right]
    \left[\begin{matrix}T;\lambda+b+c-m \\ a-m \end{matrix}\right]
    \left[\begin{matrix}b+c-m \\ b \end{matrix}\right]
  F^{(a-c-b)}\eta
\end{align*}
\color{black}
where we have used \eqnref{binomiallemma} and the calculation
\begin{align*}
T_{-1}'(&E^{(m)})T_{-1}'(F^{(b)})T_{-1}'(K^{b})
T_{-1}'(E^{(a-m)})T_{-1}'(K^{-a})
\\
&=(-1)^{b+a}v^{a(a+1) + b(b -1- 2m)}F^{(m)} E^{(b)}F^{(a-m)}
\end{align*}
By the Clebsch-Gordan decomposition \lemref{CG}
\color{black}
\begin{align*}
\rho_1&\left(\beta^{m,n}_{2r}(u^{(m)}_{i}\otimes
u^{(n)}_{j})\right)
F^{(-c)}\eta  \\
&=(-1)^jv^{m(n-j)+ni-2i(n-j)-j+r^2-\frac{(m+n)^2}{4}} 
   {\qbinom m i}{\qbinom{n}{j}} \\ \\
&\quad \times ||u^{(2r)}||^{-2}
  \frac{\left[\begin{matrix} m & n & 2r
  \\ i & {n-j} & \frac{n-m}{2}+r+i-j \end{matrix}\right]}{\left[\begin{matrix} 2r
  \\  \frac{n-m}{2}+r+i-j \end{matrix}\right]}\\   \\
&\quad \times 
(v^{2\lambda+2+i-j+\frac{n-m}{2}+4c}T^2)^{j-i+\frac{m-n}{2}}\\   \\
&\quad \times 
 \sum_{l=0}^{i-j+\frac{n-m}{2}+r} (-1)^{r-l}
  v^{l\left(r+j-i+\frac{m-n}{2}+1 \right) }
 \qbinom{i-j+\frac{n-m}{2}+c}{ i-j+\frac{n-m}{2}+r-l}
  \\
&\hskip 50pt \times \qbinom{T;\lambda+l+c}{ r}
 \qbinom{l+c}{ l}
  F^{(j-i+\frac{m-n}{2}-c)}\eta \\ 
\end{align*}
for $\max\{|\frac{m-n}{2}|,|\frac{m-n}{2}-i+j|\}\leq r\leq \frac{m+n}{2}$,
and
\begin{align*}
\rho_1&\left(T_{-1}'\beta^{m,n}_{2r}(u^{(m)}_{i}\otimes
u^{(n)}_{j})\right)
F^{(-c)}\eta  \\
&=(-1)^jv^{m(n-j)+ni-2i(n-j)-j+r^2-\frac{(m+n)^2}{4}} 
   {\qbinom m i}{\qbinom{n}{j}} \\ \\
&\quad \times ||u^{(2r)}||^{-2}
  \frac{\left[\begin{matrix} m & n & 2r
  \\ i & {n-j} & \frac{n-m}{2}+r+i-j \end{matrix}\right]}{\left[\begin{matrix} 2r
  \\  \frac{n-m}{2}+r+i-j \end{matrix}\right]}  (v^{ 2\,c
+\lambda}T)^{i-j+\frac{n-m}{2}}\\   \\
&\quad\times  \sum_{l=0}^{i-j+\frac{n-m}{2}+r}(-1)^{r-l}
  v^{  l (r +j-i+\frac{m-n}{2}+1) }
    \left[\begin{matrix}T;\lambda+c \\ l \end{matrix}\right]
    \left[\begin{matrix}T;\lambda+r+c-l \\ i-j+\frac{n-m}{2}+r-l \end{matrix}\right]
  \\
&\hskip 50pt \times    \left[\begin{matrix}r+c-l \\ r \end{matrix}\right]
  F^{(i-j+\frac{n-m}{2}-c)}\eta  
\end{align*}
\end{proof}

\begin{cor}\label{firstbetacor}  If $m=n$ and $i=j$, then we get
\begin{align*}
\rho_1&\left(\beta^{m,n}_{2r}(u^{(m)}_{i}\otimes
u^{(m)}_i)\right)
F^{(c)}\zeta  \\
&=(-1)^iv^{2i(i-m-1)+r^2} 
   {\qbinom m i}^2 \frac{\left[\begin{matrix} m & m & 2r
  \\ i & {m-i} & r \end{matrix}\right]}
   {||u^{(2r)}||^2\left[\begin{matrix} 2r
  \\  r \end{matrix}\right]} \\
&\quad \times \sum_{l=0}^{r} (-1)^{l-r}
  v^{( l-r)\left(r-1\right) }
 \\
&\hskip 10pt \times   \qbinom{T;\lambda-\rho-c+r-l}{r-l}\qbinom{l+c}{r}
\qbinom{T;\lambda-\rho-c}{l} F^{(c)}\zeta.
\end{align*}
\end{cor}

\begin{cor}\label{secondbetacor} Let $m$, $n$ be non-negative integers and
$\max\{|\frac{m-n}{2}|,|\frac{m-n}{2}-i+j|\}\leq r\leq \frac{m+n}{2}$ with $0\leq
i\leq m$ and
$0\leq j\leq n$.    Then
\begin{align*}
\rho_1\Big(\beta^{m,n}_{2r}&(u^{(m)}_i\otimes u^{(n)}_{j})
 \Big) u^{(k)}_s\\
&=(-1)^{i+\frac{n-m}{2}}v^{m(n-j)+ni-2i(n-j)-j+r^2-\frac{(m+n)^2}{4}+(j-i-
	\frac{n-m}{2})(k-2s+1)} 
   {\qbinom m i}{\qbinom{n}{j}} \\  &\quad \times ||u^{(2r)}||^{-2}
  \frac{\left[\begin{matrix} m & n & 2r
  \\ i & {n-j} & \frac{n-m}{2}+r+i-j \end{matrix}\right]}{\left[\begin{matrix} 2r
  \\  \frac{n-m}{2}+r+i-j \end{matrix}\right]}\\  &\quad \times 
 \sum_{q=0}^{i-j+\frac{n-m}{2}+r}(-1)^{ q}
  v^{- \left( i-j+\frac{n-m}{2}-r-1 \right) q}\qbinom{k+q-s}{q}\qbinom{r+s-q}{s-q}
  \\ &\hskip 100pt \times   \qbinom{k+i-j+\frac{n-m}{2}-s}{i-j+\frac{n-m}{2}+r-q}
u^{(k)}_{s-i+j-\frac{n-m}{2}}
\end{align*}
\end{cor}

\begin{proof}  The calculation is very similar to the proof of the previous proposition and so omitted.

\end{proof}
\section{Symmetry Properties of Induced Forms}

\subsection{Twisted action of $R$.}  We shall twist by an automorphism
of 
$\ru$ in the setting of $\ru$-modules. Let $\t$ be an automorphism of
$\ru$. Then for any $\ru$-module
${_R\mathcal E}$ define a new $\ru $-module ${_R\mathcal E}^\t$ with set equal to
that of ${_R\mathcal E}$ and action given by : for $e\in {_R\mathcal E}$ and $x\in\ru $
the action of $x$ on $e$ equals $\t(x)e$. For any
$\ru$-module $A$, let $A^{s\t}$ denote the module with action $\di $ on
$A^{s\t}$ defined as follows: For
$a\in A$ and $x\in\ru$, 
$$ x\di a= s\t (x)\ a\ .
$$
\smallskip

\begin{lem}
Suppose  $\phi$ is a $\brho$-invariant $R$-valued pairing of $\ru $-modules
$A$ and $B$.  Then $s\circ \phi$ is a $\brho$-invariant pairing of
$A^s$ and $B^s$.  Furthermore if either $A$ or $B$ is $\text{ad}\,F$ locally finite, then $s\circ \phi\circ L$ a
$\rho_\pi$-invariant pairing of
$A^{sT'_{-1}}$ and
$B^{sT'_{-1}}$. Also
$\phi$ and $\phi\circ L$ respectively are $\brho$-invariant
pairing of these two pairs taking values
in the $R$-module $R^s$.
.
\end{lem}

\begin{proof} By \corref{Lisomorphism} the map $L:A^{T'_{-1}}\otimes B^{T'_{-1}}\to (A\otimes B)^{T'_{-1}}$ is
an
$\ru$-module homomorphism. Now we know that $\phi\in\Hom_{\ru}(A\otimes B^{\rho_1},R)$ is module homomorphism and so
after twisting the action we get $\phi\in\Hom_{\ru}((A\otimes B^{\rho_1})^{T_{-1}'},R)$.  Composing with $s$ we
have $s\circ\phi\circ L\in \Hom_{\ru}(A^{T_{-1}'}\otimes B^{\rho_1T_{-1}'},R^s)$ is a $\bu$-module homomorphism
that is $s$-linear.  To make it $R$-linear we twist the action on $A^{T_{-1}'}$ and $B^{\rho_1T_{-1}'}$ by $s$. 
Indeed let
$\sharp$ denote the action of
$R$ twisted by
$s$. For
$r\in R,a\in A^{T_{-1}'}$ and
$b\in B^{\rho_1T_{-1}'}$, since $s^2 = 1$ , $ s\circ\phi(r\sharp a,b) = s\circ\phi(s(r)a,b)
= r\phi(a,b) = s\circ\phi(a,r\sharp b).$ 
\end{proof} 

Recall from \cite{MR96c:17022} we define a {\it cycle} ( for $A$ ) to be a pair $(A,\Psi)$
where
$A$ is a $\mathbf U$ (or $_R\mathbf U$)  module and $\Psi$ is a module
homomorphism  
\begin{equation}
\Psi : A_\pi^{sT_{-1}'} \
\ra \ A .
\end{equation} 
We say that the cycle is {\it nondegenerate} whenever $\Psi$ is an
isomorphism and that $A$ {\it admits a nondegenerate cycle} when such a
pair exists with $\Psi$ an isomorphism.

Throughout the following, suppose 
$$
\Psi:(A_\pi)^{sT_{-1}'}\to A,
$$ 
is a 
nondegenerate cycle.  For example if we return to the setting of
section one and set
$M$ equal to the Verma module with highest weight
$Tv^{-1}$: i.e.
$M=M(m)$, then $M_\pi^{sT_{-1}'}$ is isomorphic to $M$
itself.  Then we  choose maps $\Psi : M_F\ra M$. For $a\in M, x\in \ru, \
\Psi(sT_{-1}'(x)\cdot a)=x\cdot \Psi(a)$. More precisely 

{\begin{lem}\label{vermacycle} Let $m$ be a lowest weight vector in $M_\pi$ of weight
$Tv$ and $\Psi(m)$ the highest weight vector in $M$ of weight
$Tv^{-1}$ where $F^{-1}\Psi(m)\equiv m\mod \iota M$.  The map
$\Psi:(M_\pi)^{sT_{-1}'}\to M$ given by
$$
\Psi(F^{-k}m)=\frac{(-1)^kv^{-k^2}T^k}{[k]![T;-1]_{(k)}}F^{k}\Psi(m)
$$
for $k\geq 0$ is an $R$-linear $\bu$-module isomorphism.
\end{lem}}
The above can be rewritten as 
$$
\Psi(F^{(-k)}m)=F^{(k)}\Psi(m).
$$

\begin{proof}  The proof is straight forward and so omitted.
\end{proof}

Set
\begin{equation}
\bar \Psi:=\Psi \otimes sT_{1}''\circ L^{-1}:(M_\pi\otimes \mathcal E)^{sT_{-1}'}\to
M\otimes \mathcal E.
\end{equation}
Let $\iota:\mathbb P(M\otimes \mathcal E, N\otimes
\mathcal F)\to\hom_\ru(M\otimes \mathcal E\otimes (N\otimes
\mathcal F)^{\rho_1},R)$ be the canonical isomorphism with $\iota(\chi)(a\otimes b)=\chi(a,b)$.  Note that $a\otimes b\in
M\otimes
\mathcal E\otimes (N\otimes
\mathcal F)^{\rho_1}$ on the left hand side, while $(a,b)\in M\otimes \mathcal E\times N\otimes
\mathcal F$ on the right hand side.  Define $\chi\mapsto \chi^\sharp$ in $
\End(\mathbb P(M\otimes
\mathcal E, N\otimes
\mathcal F))$ by
\begin{equation}\label{twistedpairing}
\iota(\chi^\sharp)(\bar\Psi(a)\otimes \bar\Psi(b)):=
s\circ\iota(\chi_\pi)\circ L(a\otimes b)
\end{equation}
for $a\in(M\otimes \mathcal E)_\pi^{sT_{-1}'}$,
$b\in (N\otimes\mathcal  F)^{s T_{-1}'\rho_1}_\pi$ 
and $\chi\in\mathbb P(M\otimes \mathcal E,N\otimes \mathcal F)$. 
Recall that when one evaluates $L$ on the
right hand side, one has
\begin{equation*}
L(a\otimes b)=\sum_n(-1)^nv^{-n(n-1)/2}\{n\}F^{(n)}a\otimes E^{(n)}b
\end{equation*} 
where $F^{(n)}a$ is evaluated using the untwisted action and $E^{(n)}b$ is evaluated with the action twisted only by
$\rho_1$.    Note also that as a linear map
$L\in\text{End}\,((M\otimes
\mathcal E)_\pi\otimes(N\otimes
\mathcal F)_\pi^{\rho_1})$ is well defined as $F$ acts locally nilpotently on $
(M\otimes \mathcal E)_\pi$. We can view $L:(M\otimes \mathcal
E)_\pi^{sT_{-1}'}\otimes (N\otimes \mathcal
F)_\pi^{ sT_{-1}'\rho_1}\to  ((M\otimes \mathcal E)_\pi\otimes
(N\otimes \mathcal F)^{\rho_1}_\pi)^{sT_{-1}'}$ as a module isomorphism (see
\eqnref{Tandrho} and 
\corref{Lisomorphism}). Then
$s\circ \iota(\chi_\pi)\circ L:(M\otimes \mathcal
E)_\pi^{sT_{-1}'}\otimes (N\otimes \mathcal
F)_\pi^{sT_{-1}'\rho_1}\to 
R$ is a module homomorphism.  Indeed $\iota(\chi_\pi)\in\hom_\ru((M\otimes \mathcal E)_\pi\otimes (N\otimes
\mathcal F)^{\rho_1}_\pi,R)$, which implies $\iota(\chi_\pi)\in\hom_\ru(((M\otimes \mathcal E)_\pi\otimes (N\otimes
\mathcal F)^{\rho_1}_\pi)^{sT'_{-1}},R^s)$.  This gives us $\iota(\chi_\pi)\circ L\in\hom_\ru((M\otimes \mathcal
E)_\pi^{sT'_{-1}}\otimes (N\otimes
\mathcal F)^{\rho_1sT'_{-1}},R^s)$.

If we suppress the map $\iota$, then we can write \eqnref{twistedpairing} as
\begin{equation*}
\chi^\sharp(\bar\Psi(a), \bar\Psi(b)):=
s\circ \chi_\pi\circ L(a\otimes b)
\end{equation*}
were we view $a\in(M\otimes \mathcal E)_\pi^{sT_{-1}'}$ with $b\in(N\otimes \mathcal F)_\pi^{sT_{-1}'}$ on the left hand side
of the equality and  $b\in(N\otimes \mathcal F)_\pi^{sT_{-1}'\rho_1}$ on the right hand side.

More explicitly  we can show that $\chi^\sharp\in
\mathbb P(M\otimes \mathcal E, N\otimes
\mathcal F)$ by the following calculation for $x\in \ru$:
\begin{align*}
\sum\chi^\#(S(x_{(2)})\bar\Psi(a),\brho(x_{(1)})\bar\Psi(b))&=
\sum\chi^\#(\bar\Psi(T_{-1}' S(x_{(2)})a),\bar\Psi(
T_{-1}' \brho(x_{(1)})b)) \\ 
  &=\sum s\circ\chi_\pi\circ L(T_{-1}' S(x_{(2)})a\otimes 
T_{-1}' \brho(x_{(1)})b) \\ 
  &= s\circ\chi_\pi(T_{-1}' S(x)L(a\otimes b)) \\ 
  &= \mathbf e(x)s\circ\chi_\pi\circ L(a\otimes b)  \\
  &=\mathbf e(x)\chi^\#(\bar\Psi(a),\bar\Psi(b))
\end{align*}
where the third equality is from \corref{Lisomorphism}, and the fourth
equality is due to the fact that $\chi_\pi$ is $\rho$-invariant. 

\begin{lem}  Let $\phi$ be the Shapovalov form on the Verma module $M$ of highest weight $Tv^{-1}$ and $\Psi$ the cycle in \lemref{vermacycle}.  Then
$$
s\circ\phi_\pi\circ L=\phi\circ\Psi\otimes \Psi.
$$
In other words $\phi^\sharp=\phi$.
\end{lem}

\begin{proof}  Let $w$ be a fixed highest weight vector of $M$ of highest weight $Tv^{-1}$ such that $\Psi(m)=w$ with $m$ a lowest weight 
vector in $M_\pi$.
Recall $\phi^\#(\Psi(a),\Psi(b))=s\circ\phi_\pi\circ L(a\otimes b)$ for
$a,b\in M_\pi$ and that the argument above shows $\phi^\#$ is $\rho$-invariant.  Observe now
\begin{align*}
\phi^\#(w,w)&=\phi^\#(\Psi(m),\Psi(m))=s\circ\phi_\pi\circ L(m\otimes m) \\
&=s\circ\phi_\pi(m,m) \quad \text{as $m$ is a lowest weight vector} \\
&=\phi(\Psi(m),\Psi(m))=\phi(w,w) \quad  \text{by \eqnref{lift}}.
\end{align*} Since
$\phi^\#$ and
$\phi$ agree on the generator $w=\Psi(m)$ of $M$, the
$\rho$-invariance property proves that they agree everywhere.

\end{proof}

If $\phi$ is a $\rho$-invariant form on $\mathcal F\times \mathcal E$, where  $\mathcal F$ and $\mathcal E$ are two finite dimensional $\bu$-modules in
the category
$\mathcal C'$ we define in a similar manner the $\rho$-invariant form $\phi^\sharp$ on $\mathcal F\times \mathcal E$ given by $
\iota(\phi^\sharp)(T''_1(a)\otimes T''_1(b)):=
\iota(\phi)\circ L(a\otimes b)$.  If $\mathcal F=\mathcal E$ is irreducible with a highest weight vector $u^{(m)}$, then
\begin{align}\label{finaltouch}
\iota(\phi^\sharp)( u^{(m)}\otimes u^{(m)})&=\iota(\phi^\sharp)(T_1''(u^{(m)}_m)\otimes T''_1(u^{(m)}_m)) \\
&=\iota(\phi)(u^{(m)}_m\otimes
u^{(m)}_m)=\iota(\phi)(u^{(m)}\otimes u^{(m)})\notag
\end{align}
The second equality comes from the definition of $L$ and last equality is due to \eqnref{normalizedforms}.  Thus
$\phi^\sharp=\phi$.  This implies that if
$\mathcal F=\oplus_i\mathcal F_{n_i}$ with
$n_i$ distinct and $\mathcal E=\oplus_j\mathcal F_{m_j}$ with
$m_j$ distinct nonnegative integers, then $\phi=\sum_{i,j,n_i=m_j}\phi_i$ where $\phi_i$ is a $\rho$-invariant nondegenerate form on $\mathcal F_{n_i}$ and
zero on the other summands $\mathcal F_{n_k}$, $k\neq i$, and we still have $\phi^\sharp=\phi$.  

For four finite dimensional $X$-admissible $\bu$-modules $\mathcal E$, $\mathcal F$, $\mathcal M$ and $\mathcal N$ and invariant form $\phi$ on $\mathcal
M\times \mathcal N$ we can define for each
$\bu$-module homomorphism $\b:{\mathcal E}\otimes \mathcal
F^{\rho_1}\ra\bu$ the {\it induced form}
$\phi_{\beta}$ by the formula, for
$e\in\mathcal E,f\in \mathcal F, m\in\mathcal M,n\in \mathcal N$,
\begin{equation}
\phi_{\beta}(m\otimes e,n\otimes f)=\phi( m,\rho_1(\beta(e\otimes f))n).
\end{equation}
This is similar to the definition \eqnref{inducedform}.
 
\begin{lem}\label{tediouslemma}
For $m$, $n$, $\sigma$, $l$ and $r$ non-negative integers with $0\leq \sigma\leq n$ and $|m-n|\leq 2r\leq m+n$, one has
\begin{align}
&v^{\frac{2n+n^2-(2\,m+m^2)}{4}+n- 2\,\sigma +r(r+1)}\\ 
&\quad \times\sum_{k=0}^{ n-\sigma}  (-1)^{k}v^{-\frac{k(k+2l+3)}{2}}(v-v^{-1})^{k} \notag \\
&\hskip 50pt \times \qbinom{r-k}{l}\frac{[r+k]![\sigma+k]_{(k)}}{[k]!}
     {\qbinom{n}{\sigma+k}}  
      \left[\begin{matrix} m & n & 2r
     \\ \frac{m-n}{2}+\sigma & {n-k-\sigma} & r-k \end{matrix}\right]\notag \\ \notag \\
&=(-1)^{\frac{m+n}{2}+r}\sum_{k=0}^{n-\sigma}v^{\frac{k\left(3+k-2l+4\sigma-2n+2r\right)}{2}}
   (v-v^{-1})^{k}  \notag\\
&\hskip 70pt \times \qbinom{r-k}{l-k}\frac{[r+k]![\sigma+k]_{(k)}}{[k]!}
     {\qbinom{n}{\sigma+k}}
   \left[\begin{matrix} m & n & 2r
  \\ \frac{m+n}{2}-\sigma & {\sigma+k} & r+k \end{matrix}\right].\notag
\end{align}

\end{lem}

\begin{proof}  The result follows from a rather tedious calculation using $\phi=(\enspace,\enspace)$, the normalized nondegenerate form on $\mathcal
F_{a}$, and 
\corref{secondbetacor}.
\end{proof}

\subsection{} In this section we suppose $\mathcal F_m$
and
$\mathcal F_n$ are the $X$-admissible finite dimensional $\bu$-modules
given in \secref{CG}.  For any homomorphism
$\beta:_R\mathcal F_m\otimes _R\mathcal F_n^{\rho_1}\to {_RF(\bu)}$ which has the
form
\begin{equation}\label{factorization}
\beta=\sum_{m,n,k}r^{m,n}_{k}\beta^{m,n}_{2k}
\end{equation}
where $r^{m,n}_{k}\in R$, observe that
$s\circ \beta$ has the same form.

\begin{thm}\label{firstinvariance}
Let $M$ be the Verma module of highest weight $Tv^{-1}$ (so that $\lambda=0$) and assume that $\beta:_R\mathcal F_m\otimes _R\mathcal F_n^{\rho_1}\to
{_RF(\bu)}$  has the form \eqnref{factorization}.
If
$\phi$ is a $\ru$-invariant pairing on $M$ satisfying $s\circ
\phi_{\pi}\circ L=\phi\circ (\Psi\otimes
\Psi)$, then 
\begin{equation}\label{hardprop}
\chi _{\b,\phi}^\sharp = \chi_{s\beta ,\phi}\ .
\end{equation} 
\end{thm}

\begin{proof}  Since $\chi _{\b,\phi}^\sharp$ is $R$-linear in $\beta$ we can
reduce to the case that $\beta=a\beta^{m,n}_{2r}$ for some $a\in R$.   
Let us choose basis $\{F^{(-k)}m_{0,1}|\,k\in \mathbb N\}$, $\{u_i^{(m)}|\,0\leq i\leq
m\}$,
$\{u_i^{(n)}|\,0\leq i\leq n\}$, respectively for $M=N$,
$\mathcal F_m$ and $\mathcal F_n$, where $m_{0,1} $ is a lowest weight vector of weight $Tv$ that generates $M_\pi$, and $u^{(m)}$ and $u^{(n)}$
are highest weight vectors of weights $v^m$ and $v^n$.  

Now the element
$m_{0,1}\otimes u^{(m)}$ generates $(M\otimes \mathcal F_m)_\pi^{sT'_{-1}}$
(see \lemref{Mackey}), so we can reduce the proof of the identity further to
showing that
\begin{equation}\label{firstreduction}
\chi _{\b,\phi}^\sharp(\bar\Psi(m_{0,1}\otimes u^{(m)}),\bar\Psi(F^{(-k)}\zeta \otimes
u_i^{(n)}))=
\chi_{s\beta ,\phi}(\bar\Psi(m_{0,1}\otimes u^{(m)}),\bar\Psi(F^{(-k)}\zeta \otimes
u_i^{(n)}))
\end{equation}
for all $k\in \mathbb N$ and $0\leq j\leq n$.  Indeed since both sides are
invariant, if $w\in (M\otimes \mathcal F_m)_\pi^{sT'_{-1}}$, then there
exists $u\in \ru$ such that $w=u\bar\Psi(m_{0,1}\otimes u^{(m)})$ and hence
\begin{align*}
\chi _{\b,\phi}^\sharp(w,\bar\Psi(F^{(-k)}m_{0,1} \otimes u^{(n)}_j))
&=\chi_{\b,\phi}^\sharp(u\bar\Psi(m_{0,1}\otimes u^{(m)}),F^{(-k)}m_{0,1}
\otimes u^{(n)}_j) \\  
&=\chi_{\b,\phi}^\sharp(\bar\Psi(m_{0,1}\otimes u^{(m)}),\,
   \rho(u)\bar\Psi(F^{(-k)}m_{0,1} \otimes u^{(n)}_j)) \\ 
&= \chi_{s\beta ,\phi}(\bar\Psi(m_{0,1}\otimes u^{(m)}),
 \rho(u)\bar\Psi(F^{(-k)}m_{0,1} \otimes u^{(n)}_j)) \\ 
&= \chi_{s\beta ,\phi}(u\bar\Psi(m_{0,1}\otimes u^{(m)}),
   \bar\Psi(F^{(-k)}m_{0,1} \otimes u^{(n)}_j))  \\
&= \chi_{s\beta ,\phi}(w,
   \bar\Psi(F^{(-k)}m_{0,1} \otimes u^{(n)}_j)).
\end{align*}
Here the third equality would be true as $\rho(u)\bar\Psi(F^{(-k)}m_{0,1} \otimes u^{(n)}_j)$ is a linear combination of other $\bar\Psi(F^{(-l)}m_{0,1}
\otimes u^{(n)}_s)$.
 Recall that $\mathcal K$ is the field of fractions of $R$.  By \cite[Lemma 4.2]{MR96c:17022,MR96c:17037}, to prove \eqnref{firstreduction}, we need to
show
\begin{align}\label{secondreduction}
_{\mathcal K}\chi _{\b,\phi}^\sharp(\bar\Psi(m_{0,1}\otimes u^{(m)}),\bar\Psi&(F^{(-k)}m_{0,1} \otimes
u_i^{(n)}))= \\
&{_{\mathcal K}\chi_{s\beta ,\phi}}(\bar\Psi(m_{0,1}\otimes u^{(m)}),\bar\Psi(F^{(-k)}m_{0,1} \otimes
u_i^{(n)}))\notag
\end{align}
Now $_{\mathcal K}M_\pi\otimes \mathcal E_n$ is a direct sum of lowest weight irreducible Verma modules with lowest weight vector $m_{\epsilon
r',\epsilon r'+1}$ of weight $Tv^{\epsilon r'+1}$, where
$\epsilon=\pm 1$ and $r'=n-2p$, $0\leq p\leq n$.  From this fact we get that $F^{(-k)}m_{0,1} \otimes u^{(n)}_j$ is a $\mathcal K$-linear combination
of 
$E^{(i)}m_{\epsilon r',\epsilon r'+1}$.
Hence we need to show
\begin{align*}
_{\mathcal K}\chi _{\b,\phi}^\sharp(\bar\Psi(m_{0,1}\otimes u^{(m)}),\bar\Psi&(E^{(i)}m_{\epsilon r',\epsilon r'+1}))= \\
&{_{\mathcal K}\chi_{s\beta ,\phi}}(\bar\Psi(m_{0,1}\otimes u^{(m)}),\bar\Psi(E^{(i)}m_{\epsilon r',\epsilon r'+1}))
\end{align*}
for all $r'$, $\epsilon$ and $i$ with $\epsilon r'+2i=m$.  Using the invariance of the two forms and the fact that $m_{0,1}$ is a lowest weight vector,
proving the above identity is equivalent to verifying
\begin{equation*}
_{\mathcal K}\chi _{\b,\phi}^\sharp(\bar\Psi(m_{0,1}\otimes u^{(m)}_i),\bar\Psi(m_{\epsilon r',\epsilon r'+1}))=
{_{\mathcal K}\chi_{s\beta ,\phi}}(\bar\Psi(m_{0,1}\otimes u^{(m)}_i),\Psi(m_{\epsilon r',\epsilon r'+1})).
\end{equation*}

Recall the definition of $\sigma(\epsilon,n,p)$ and $a_{\epsilon, n, p}(k)$ from \eqnref{sigma}.  Then
\begin{align*}
\chi_{\beta,\phi}^\sharp(&\bar\Psi(m_{0,1} \otimes u^{(m)}_i),
  \bar\Psi(m_{\epsilon r',\epsilon r'+1}))=
s\circ\chi_{\beta,\phi_\pi}
   \circ L(m_{0,1} \otimes u^{(m)}_i,m_{\epsilon r',\epsilon r'+1})\\
&= \sum_q(-1)^qv^{-\frac{q(q-1)}{2}}\{q\}s\circ\chi_{\beta,\phi_\pi}(
   F^{(q)}(m_{0,1} \otimes u^{(m)}_i),
  E^{(q)}*(m_{\epsilon r',\epsilon r'+1})) \\
&= s\circ\chi_{\beta,\phi_\pi}(m_{0,1} \otimes u^{(m)}_i,
   \mathcal L
   (m_{\epsilon r',\epsilon r'+1})) \\
&= s\circ\chi_{\beta,\phi_\pi}(m_{0,1} \otimes u^{(m)}_i,m_{\epsilon r',\epsilon r'+1}) \\
&= s\circ\phi_\pi(m_{0,1} , \sum_k a_{\epsilon, n,p}(k)(\epsilon)\rho_1(\beta(u^{(m)}_i\otimes
  u^{(n)}_{k+\sigma}))F^{(-k)}m_{0,1} ) \\
&= \phi(\Psi(\mathcal L^{-1}m_{0,1} ), 
   \Psi(\sum_k a_{\epsilon, n,p}(k)\rho_1(\beta(u^{(m)}_i\otimes
  u^{(n)}_{k+\sigma}))F^{(-k)}m_{0,1} )) \\ 	
&= \phi\left(\Psi(m_{0,1} ), 
   \Psi\left(\sum_k a_{\epsilon,n,p}(k)\rho_1(\beta(u^{(m)}_i\otimes
  u^{(n)}_{k+\sigma}))F^{(-k)}m_{0,1}
  \right)\right)	
\end{align*}
 
The seventh equality is due to the fact that $\mathcal L^{-1}m_{0,1}=m_{0,1}$.   We may assume that $\sigma-i+\frac{m-n}{2}=0$, 
as distinct root spaces are orthogonal. On the other hand 
\begin{align*}
\chi_{s\beta,\phi}&(\bar\Psi(m_{0,1} \otimes u^{(m)}_i), 
  \bar\Psi (m_{\epsilon r,\epsilon r+1}))  \\
&= v^{-2(n-\sigma)(\sigma+1)}\sum_{s=0}^{n-\sigma}v^{s(1-\epsilon (n-2p))}a_{\epsilon, n, p}(k) \\
&\hskip 30pt \times \phi\left(\Psi(m_{0,1} ),
\Psi\left(\left(\rho_1\circ
T'_{-1}\circ\beta
  \right)\left(T''_1(u^{(m)}_i)\otimes T''_1(u^{(n)}_{\sigma+k}\right)\right) F^{(-k)}m_{0,1})\right)
\end{align*} 
We used the reduction to $s\beta =\beta$ and the fact that $\rho_1\circ T'_{-1}=T'_{-1}\circ \rho_1$. Consequently the two sides of \eqnref{hardprop} are
equal provided
\begin{align}
\sum_{k=0}^{n-\sigma} &a_{\epsilon,n,p}(k)\rho_1(\beta(u^{(m)}_i\otimes
  u^{(n)}_{k+\sigma}))F^{(-k)}m_{0,1} \\
&= v^{-2(n-\sigma)(\sigma+1)}\sum_{k=0}^{n-\sigma}v^{s(1-\epsilon (n-2p))}a_{\epsilon, n, p}(k)\notag \\
&\hskip 50pt\times\left(\rho_1\circ
T'_{-1}\circ\beta
  \right)\left(T''_1(u^{(m)}_i)\otimes T''_1(u^{(n)}_{\sigma+k})\right) F^{(-k)}m_{0,1}.\notag
\end{align}

We will now expand out the two sides and eventually show that they are equal.
We begin with a calculation for the left hand side. 
\begin{align*}
\rho_1\Big(\beta^{m,n}_{2r}(&u^{(m)}_i\otimes u^{(n)}_{k+\sigma})
 \Big) F^{(-k)}m_{0,1}\\
&=(-1)^{\sigma}v^{k(2k+ 2i - m )  +(m- i)n - \frac{(m+n)^2}{4} +  (2i-1- m)\sigma + r(1 + 2r)} 
   {\qbinom m i}{\qbinom{n}{k+\sigma}} \\ 
&\quad \times ||u^{(2r)}||^{-2}
  \left[\begin{matrix} m & n & 2r
  \\ i & {n-k-\sigma} &r-k \end{matrix}\right]\frac{[r-k]![r+k]!}{ [2r]!}T^{2k}
  \qbinom{T;r}{ r}\qbinom{r}{k}m_{0,1}
\end{align*}
This implies
\begin{align*}
&\sum_k a_{\epsilon,n,p}\rho_1(\beta(u^{(m)}_i\otimes
  u^{(n)}_{k+\sigma}))F^{(-k)}m_{0,1}
   \\
&=(-1)^{\sigma}\frac{v^{(2i-1- m)\sigma + r(1 + 2r)+(m- i)n - \frac{(m+n)^2}{4}}}{||u^{(2r)}||^2}
{\qbinom m i} 
[T;r]_{(r)}\\ 
&\quad \times\sum_{k=0}^{\min\{n-\sigma,r\}}  (-1)^k 
  v^{-k} 
   {\qbinom{n}{k+\sigma}}
  \frac{[\sigma+k]_{(k)}}{
[T;k]_{(k)}}\frac{[r+k]!}{ [k]![2r]!}\left[\begin{matrix} m & n & 2r
  \\ i & {n-k-\sigma} & r-k \end{matrix}\right]T^{k}m_{0,1}
\end{align*}

by \lemref{highestlowestwtvectors} and \corref{betaonabasis}. 
The right hand side is 
\begin{align*}
&v^{-2(n-\sigma)(\sigma+1)}\sum_{k=0}^{n-\sigma}(-1)^k\frac{v^{k(1-n+2\sigma)}[\sigma+k]_{(k)}}{T^k[T,k]_{(k)}} \\
&\hskip 100pt\times \left(\rho_1\circ T'_{-1}\circ\beta
  \right)\left(T''_1(u^{(m)}_i)\otimes T''_1(u^{(n)}_{\sigma+k})\right) F^{(-k)}m_{0,1}  \\ 
&\quad\times  \sum_{l=0}^{k+r}(-1)^{r-l}
  v^{  l (r -k+1) }
    \left[\begin{matrix}T;k \\ l \end{matrix}\right]
    \left[\begin{matrix}T;r+k-l \\ r+k-l \end{matrix}\right]
    \left[\begin{matrix}r+k-l \\ r \end{matrix}\right]m_{0,1}  \\ 
&=(-1)^{m-i+r}\frac{v^{i(i-1)+ m - 2\,n + m\,n + \sigma( \sigma+2 - m - n ) + r^2 - \frac{(m+n)^2}{4}}}{||u^{(2r)}||^2}{\qbinom{m}{i}}
\left[T; r \right]_{(r)}\\   
&\quad \times\sum_{k}(-1)^k\frac{[\sigma+k]_{(k)}
  v^{k\left(1 + 2\,i - m \right)}}{T^k
[T;k]_{(k)}}
   {\qbinom{n}{\sigma+k}} 
  \frac{[r+k]!}{ [k]![2r]!}\left[\begin{matrix} m & n & 2r
  \\ m-i & {\sigma+k} & r+k \end{matrix}\right]
    m_{0,1}  \\
\end{align*}
which follows from the Chu-Vandermonde formula, \lemref{chuvandermonde}.

Thus to prove they are equal we need for $2\sigma-2i+m-n=0$ and $|m-n|\leq 2r\leq m+n$, that
\begin{align*}
&v^{(2\sigma-n-1)\sigma + r(1 + r)+ \left(\frac{m+n}{2}-\sigma\right)n }\\ 
&\quad \times\sum_{k=0}^{n-\sigma}  (-1)^k 
  v^{-k} 
  \frac{[\sigma+k]_{(k)}}{[r]!}
     {\qbinom{n}{\sigma+k}}\qbinom{r}{k}\frac{[r-k]![r+k]!}{[2r]!}
\\
&\hskip 50pt\times \left[\begin{matrix} m & n & 2r
  \\ \frac{m-n}{2}+\sigma & {n-k-\sigma} & r-k \end{matrix}\right][T;r]_{(r-k)}T^{k}\\ 
&=(-1)^{\frac{m+n}{2}+r}\sum_{k=0}^{\min\{r, n-\sigma\}}(-1)^kv^{k\left(1 + 2\sigma-n \right)}\frac{[\sigma+k]_{(k)}}{[r]!}
     {\qbinom{n}{\sigma+k}}\qbinom{r}{k}\frac{[r-k]![r+k]!}{[2r]!}	\\
&\hskip 50pt \times\left[\begin{matrix} m & n & 2r
  \\ \frac{m+n}{2}-\sigma & {\sigma+k} & r+k \end{matrix}\right][T;r]_{(r-k)}T^{-k}.
\end{align*}
In order to prove that these are equal we need to expand in powers of $T$ and then show that the coefficients are
equal:
\begin{align*}
&v^{\frac{2\,n+n^2-(2m+m^2)}{4}+ n-2\sigma + \frac{3r(r+1)}{2}}\\ 
&\quad \times\sum_{l\geq 0}(-1)^lv^{-l(r+1)}\sum_{k=0}^{n-\sigma}  (-1)^{k} 
 \frac{[r+k]!}{[n-\sigma-k]![k]!}\left[\begin{matrix} m & n & 2r
  \\ \frac{m-n}{2}+\sigma & {n-k-\sigma} & r-k \end{matrix}\right] \\ 
&\hskip 100pt \times\frac{v^{-\frac{k(k+2l+3)}{2}}}{ (v-v^{-1})^{r-k}} \qbinom{r-k}{l}T^{r-2l}  \\ 
&=(-1)^{\frac{m+n}{2}+r}v^{\frac{r(r+1)}{2}}\sum_{l\geq 0}(-1)^lv^{-l(r+1)}\sum_{k=0}^{ n-\sigma}
 \frac{v^{\frac{k\left(3+k-2l +
4\sigma-2n+2r\right)}{2}}}{(v-v^{-1})^{r-k}} \frac{[r+k]!}{[n-\sigma-k]![k]!} \\
&\hskip 150pt \times
  \left[\begin{matrix} m & n & 2r
  \\ \frac{m+n}{2}-\sigma & {\sigma+k} & r+k \end{matrix}\right]\qbinom{r-k}{l-k}T^{r-2l} \\
\end{align*}
This reduces to the follwing identity between Clebsch-Gordan coefficients that must hold for all $m$, $n$ with $m+n$ even, $0\leq l\leq r$, $|m-n|\leq 2r\leq
m+n$,
$0\leq \sigma\leq n$:
\begin{align}
&v^{\frac{2n+n^2-(2\,m+m^2)}{4}+n- 2\,\sigma +r(r+1)}\label{byinduction}\\ 
&\quad \times\sum_{k=0}^{ n-\sigma}  (-1)^{k}v^{-\frac{k(k+2l+3)}{2}}(v-v^{-1})^{k}\frac{[\sigma+k]_{(k)}[r+k]!}{[k]!}\notag \\
&\hskip 100pt \times     {\qbinom{n}{\sigma+k}}\qbinom{r-k}{l} 
      \left[\begin{matrix} m & n & 2r
     \\ \frac{m-n}{2}+\sigma & {n-k-\sigma} & r-k \end{matrix}\right]\notag \\ \notag \\
&=(-1)^{\frac{m+n}{2}+r}\sum_{k=0}^{n-\sigma}v^{\frac{k\left(3+k-2l+4\sigma-2n+2r\right)}{2}}
   (v-v^{-1})^{k}\frac{[\sigma+k]_{(k)}[r+k]!}{[k]!}  \notag\\ 
&\hskip 100pt \times     {\qbinom{n}{\sigma+k}}
  \qbinom{r-k}{l-k} \left[\begin{matrix} m & n & 2r
  \\ \frac{m+n}{2}-\sigma & {\sigma+k} & r+k \end{matrix}\right].\notag
\end{align}

The above identity follows from \lemref{tediouslemma}.

\end{proof}

\subsection{}\label{decompositionsection} Let us return to the setting of section one.  Fix a finite dimensional  $X$-admissible 
 $\bu$-module $\mathcal F$  with highest weight $v^n$ and let $M$ be
the Verma module of highest weight $Tv^{-1}$.  Recall from \cite[{\S 2}]{MR96c:17022} the
modules $P(m+\l):=P_{m+\la}$. Then we have the decomposition 
$M\otimes_R\mathcal F =\sum_iP(m+i)$ where the sum is over the nonnegative
weights of $\mathcal F$ and  by convention we set $P(m)=M(m)$. Set
$P_i=P(m+i)$ and following the notation of \cite[3.6]{MR96c:17037} let $\mathbb Z_i$ equal the
set of integers with the opposite parity to $i$. For $j\in \mathbb Z_\i$ , set
$z^i_j = \vj + \vsj$. Then for $i\in \mathbb N^*$, the set $\{ [T;0] \vj :j\in
\mathbb Z_i \} \cup \{ z^i_j : j\in \mathbb Z_\i \}$ is an $R$ basis for the 
localization $P_{i,F}$. Also the action of $\rui$ is given by the formulas
in \lemref{locstructure} as well  as the formulas : for all indices
$j\in
\mathbb Z_\l$,
\begin{gather*}\label{Prelations}
K_\mu z_j^i = T\, v^{\angles{\mu, s\l-\rho-l_ji'}}z_j^i, \\ F
z_j^i = z_{j-2}^i, \\ E z_j = [l_j][T;-k_j] z_{j+2}^i +[k_j-l_j][T;0] w_{\l ,j+2}.
\end{gather*}

Fix a positive weight $v^r$ of $\mathcal F$ and let $P=P_r$. Set 
$\mathfrak L$ equal to
the
$m-r$-th weight space of $P$. Then $\mathfrak L$ is a free rank two $R$-module with
basis $\{z^r_{-r-1},[T;0] w_{r,-r-1}\}$. Define an $s$-linear map
$\G$ on $\mathfrak L$ and constants $a_{\pm r}$ by the formula: 
\begin{equation}\label{ar}
\G ([T;0] w_{\e r,-r-1}) = \bar\Psi([T;0] w_{\e r,r+1}) = a_{\e r}\ [T;0]
w_{-\e r,-r-1}.
\end{equation}

\def\ov{\overline}

This s-linear map $\G$ is the mechanism by which we analyze the symmetries
which arise through the exchange of
$\mathfrak L\cap ([T;0]\cdot M(m+r))$ and $\mathfrak L\cap ([T;0]\cdot M(m-r))$. This is a
fundamental calculation for all which follows. Set $\overline{\G}=[r]!\
\G$.

\begin{lem}\label{exchange} Let $\epsilon=\pm 1$.  For $a,b\in \mathfrak L\cap M(m+\epsilon r)$, we have: 
$$
\chi^\sharp(\G a,\G b)=\frac{1}{ T^r [r]![T^{-\epsilon};r]_{(r)}}\ s\chi(a,b)\quad
and \quad
\chi^\sharp(\ov\G a,\ov\G b)=u_\epsilon\ s\chi(a,b),
$$  where $u_\epsilon$ is a unit and $u_\epsilon\equiv 1 \mod (T-1)$. 

\end{lem}
\begin{proof} Using the defining identity for localization
\cite[4.3.2]{MR96c:17022}, we obtain, for any invariant form
$\chi$ on
$P$ and
$\epsilon=\pm 1$,
\begin{align*}
a_{\epsilon r}^2\ \chi^\sharp ([T;0] w_{-\epsilon r,-r-1},&[T;0]
w_{-\epsilon r,-r-1})
        = \chi^\sharp (\G [T;0] w_{\epsilon r,-r-1},\G [T;0]
w_{\epsilon r,-r-1}) \\
        &=s\chi_{_F}([T;0] w_{\epsilon r,r+1}, [T;0]
w_{\epsilon r,r+1})  \\
        &=\frac{1}{T^r [r]![T^{-\epsilon};r]_{(r)}} s\chi([T;0]
w_{\epsilon r,-r-1}, [T;0] w_{\epsilon r,-r-1}) \ .
\end{align*}  This implies the lemma.
\end{proof} 

\subsection{} Now we turn to the delicate calculation of the constant
$a_{\pm r}$.

\begin{lem}\label{delicate} We may choose a basis for $P_r$ satisfying the relations in
\eqnref{locstructure} and \eqnref{Prelations}, 
dependent only on the cycle $\bar\Psi$, and such that the
constants $a_{\pm r}$ are uniquely determined by the three relations:
$$ 
a_{-r}=s\ a_r\ , \quad a_r^2=\frac{1}{
[r]![T^{-1};r]_{(r)}}\quad and
\quad a_r\equiv \frac{ (-1)^{r+1}}{ [r]!}
\mod\ T-1\ .
$$
\end{lem}
\begin{cor}\label{matrixcor} For $\epsilon=\pm 1$,
$\overline\G ([T;0] w_{\epsilon r,-r-1}] = b_{r,\epsilon} [T;0] w_{-\epsilon r,-r-1}\ $,  where
$b_{r,\epsilon}$ is the unit determined by conditions: 
$$ 
b_{r,\epsilon}^2=\frac{[r]!}{ [T^{-\epsilon};r]_{(r)}}\quad and
\quad b_{r,\epsilon}\equiv -1-\epsilon\a(r)(T-1)\ \mod\ (T-1)^2\ 
$$ 
where $\displaystyle{\a(r)=\frac{r}{2}-\{1\}^{-1}\sum_{s=1}^r
\frac{v^s}{[s]}}$ (see \eqnref{matrixentry}).  Moreover $\ov\G$ induces a $k(v)$-linear map on
$\mathfrak L/
\ (T-1)\cdot
\mathfrak L$ given by the matrix
\begin{equation}
\begin{pmatrix}
 1&-\a(r)\{1\} \\  0&1
\end{pmatrix}.
\end{equation}
  Moreover, if
$x_\epsilon\in M(m+\epsilon r)$ and $[T;0]\cdot x_\epsilon$ is an $R$-basis vector for $\mathfrak L\cap
M(m+\epsilon r)$, then $\{[T;0]\cdot x_\epsilon,x_\epsilon +\overline \G x_ \epsilon\}$ is an $R$-basis
for $\mathfrak L$ and $x_\epsilon+\overline \G x_\epsilon$ generates the
$\ru $-submodule $P_r$.
\end{cor}
\begin{proof} The first identity is a direct consequence of the
lemma. As for the second set
$z=z^r_{-r-1}$ and $w_\pm=w_{\pm r,-r-1}$ and let $\pi$ denote the
projection $\pi:\mathfrak L\ra \mathfrak L/\ (T-1)\cdot \mathfrak L$.  Then $\mathfrak L$ has $R$-basis $\{[T;0]
w_+,z\}$ and since $z=w_++w_-$, we obtain $\pi([T;0] w_+)= -\pi([T;0] w_-)$.
So $\overline \G([T;0] w_+)=b_{r,1}[T;0] w_-\equiv -b_{r,1}[T;0] w_+\equiv [T;0]
w_+$ mod 
$(T-1)\cdot \mathfrak L$. This gives the first column in the matrix.

From the $s$-linearity of $\ov\G$ and the congruence for the unit 
$b_{r,\epsilon}$ we obtain: 
\begin{align*}
\ov\G z&=-b_{r,1}w_--b_{r,-1}w_+\equiv z +\a(r)((T-1) w_--(T-1) w_+) \\
&\equiv z-2\a(r)(T-1) w_+
\equiv z- \a(r)\{1\}[T;0] w_+ \mod (T-1)\cdot \mathfrak L
\end{align*} 

To prove the identity for $\alpha(r)$ just solve $\displaystyle{b_{r,\epsilon}^2=\frac{[r]!}{ [T^{-\epsilon};r]_{(r)}}}$ using \eqnref{matrixentry}.
Finally a short calculation shows that the transition matrix from the basis
$\{z,[T;0]w_\epsilon\}$ to 
$\{w_\epsilon +\overline \G w_\epsilon,[T;0]w_\epsilon\}$ is 
\begin{equation}\label{transistionmatrix}
 \begin{pmatrix}
1+[T;0]\overline u_\epsilon &0\\-\overline u_\epsilon&1
\end{pmatrix} 
\end{equation}
where $\overline u_\epsilon$ is the unit
$$
(-b_{r,\epsilon}-1)/[T;0]\equiv \epsilon\frac{\a(r)T(v-v^{-1})}{ 1+T}\equiv \frac{\epsilon}{
2}\a(r)\{1\}\mod (T-1).
$$  If $[T;0]\cdot x$ is
a basis vector of $\mathfrak L\cap M(m+\epsilon r)$ then
$x=uw_\epsilon$ for some unit $u$ and the transition matrix takes the
form as in \eqnref{transistionmatrix},  with the unit $\overline u_\epsilon$ replaced by
the element of $R$ equal to
$(-(su)b_{r,\epsilon}-u)/[T;0]$.  The determinant of this matrix is a
unit which implies we have a basis.  This completes the proof of the
corollary.
\end{proof}

\begin{proof} We begin with any basis satisfying \lemref{locstructure} and \eqnref{Prelations}
and let $a_{\pm r}$ be given by \eqnref{ar}. Since the vector $z^r_{r+1}$
projects to a cyclic vector of its generalized eigenspace in 
$M_\pi\otimes_R \mathcal F$, we find that
$\Psi(z^r_{r+1})=-a_rw_{-r,-r-1}-a_{-r}w_{r,-r-1}$ is a cyclic vector of
$P$. 
Thus both
$a_{\pm r}$ are units and $a_r\equiv a_{-r}\mod T-1$.
Any automorphism $\kappa$
is given by multiples of the identity $\beta\cdot 1$ on $M(m+r)$ and
$\gamma\cdot 1$ on $M(m-r)$ with both $\beta$ and $\gamma$  units and
$\beta\equiv \gamma\mod T-1$. Set $w_\epsilon=w_{\epsilon r,-r-1}$.  Then $\G(\kappa
w_+)=\G(\beta w_+)=s\beta\G (w_+)= -\frac{s\beta}{ \gamma}a_r\kappa
w_-$ and similarly,
$\G(\kappa w_-)= -\frac{s\gamma}{ \beta }a_{-r}\kappa w_+$. From the
preceding remark
$\frac{sa_r}{ a_{-r}}\equiv 1
\mod T-1$ and so we set $s\gamma$ equal to the square root of this quotient
which is also congruent to $1\mod T-1$. Put $\beta=1$. If $\kappa$ is the
automorphism of $P$ corresponding to this choice of $\beta$ and $\gamma$,
then  applying $\kappa$ to the original basis of $P$ gives a basis for which
the constants $a_{\pm r}$ satisfy the  first identity of the lemma $
a_{-r}=s\ a_r$.

Fix a form $\phi_M$ on $M$ with $\phi_M=\phi_M^\sharp$ and let $\phi_{\mathcal F}$ be an invariant form on 
$\mathcal F$. Set
$\phi=\phi_M\otimes \phi_{\mathcal F}$. Using the invariance (\propref{invariance}) for $\phi_\pi$
we now check that
 $\phi=\phi^\sharp$.  As in the proof of \thmref{firstinvariance} we need only check on certain basis vectors:
\begin{align*}
\phi^\sharp(&\bar\Psi(m_{0,1}\otimes u^{(n)}),\bar \Psi(F^{(-i)}m_{0,1}\otimes u^{(n)}_i))= s\phi_\pi\circ L(m_{0,1}\otimes u^{(n)}\otimes 
F^{(-i)}m_{0,1}\otimes u^{(n)}_i) \\
&=\sum_pv^{-p(p-1)/2}\{p\} \\ 
&\qquad \sum_{p=r'+r''}v^{p-r'r''} s[\phi_{M,\pi}(m_{0,1},F^{(r')}F^{(-i)}m_{0,1})
  \phi_{\mathcal F}( u^{(n)}_p,K^{-r'}F^{(r'')}u^{(n)}_i) ] \\
&\hskip 100pt \text{by \cite[Theorem 5.1]{MR96c:17022}}\\ 
&=v^{2\,i^2 - i\,n}\sum_pv^{\frac{-p(p-3)}{2} + p( i - n)+p^2}
  \{p\}\qbinom{p}{p-i}\qbinom{n}{p}T^{-i}\qbinom{T^{-1};i}{i}s[\phi_{M,\pi}(m_{0,1},m_{0,1}) ] \\ 
&\hskip 100pt\text{ by \lemref{binomiallemma}  and \eqnref{normalizedforms}} \\ 
&=v^{2\,i^2 - i\,n}
\sum_pv^{\frac{p(p+3)}{2}+p(i-n)}\{p\}\qbinom{n-i}{p-i}T^{-i}\qbinom{n}{i}\qbinom{T^{-1};i}{i}\phi_{M}(\Psi(m_{0,1}),\Psi(m_{0,1})).
\end{align*}
On the other hand
\begin{align*}
&\phi(\bar\Psi(m_{0,1}\otimes u^{(n)}),\bar\Psi(F^{(-i)}m_{0,1}\otimes u^{(n)}_i)) \\
&=v^{\frac{i(i-1)}{2}+i^2+2n}\{i\}T^{-i}\qbinom{n}{i}\qbinom{T^{-1};i}{i}\phi_M(\Psi(m_{0,1}),\Psi(m_{0,1}))\phi_{\mathcal F}
  (u^{(n)},u^{(n)}).
\end{align*}
Now if we use \corref{link} with $p=0$ and replace $k$ by $p$ and $s$ by $i$ we get
\begin{equation*}
v^{i(1-r)-2n}\{i\}=\sum_{p=i}^{n}(-1)^{p+i}v^{p(r-2p-1)+\frac{(p - i)(3p + i-1) }{2}}\{p\}\qbinom{n-i}{p-i}.
\end{equation*}
Using the automorphism of $k(v)$ induced by $v\mapsto v^{-1}$ we get
\begin{align*}
v^{2n-i(1-r)}\{i\}
&=v^{\frac{i^2-i}{2}}\sum_{p=i}^{n}v^{ \frac{p^2+3p}{2} +  p\,\left(i - r\right)}\{p\}\qbinom{n-i}{p-i}.
\end{align*}
Taking $r=n$ then 
\begin{equation*}
v^{2i^2-i\,n }\sum_pv^{\frac{p(p+3)}{2}+p(i-n)}\{p\}\qbinom{n-i}{p-i}
=v^{2n+i^2+\frac{i^2-i}{2}}\{i\}.
\end{equation*}
This completes the proof that $\phi=\phi^\sharp$.

Now from \eqnref{ar}, 
\begin{align*}
a_{\epsilon r}^2\ \phi^\sharp(w_{-\epsilon r,-r-1},w_{-\epsilon r,-r-1}) & =
s\phi_\pi(w_{\epsilon r,r+1}, w_{\epsilon r,r+1}) \\
&=\frac{1}{ [r]!T^r\ [T^{-\epsilon};r]_{(r)}}s\phi(w_{\epsilon r,-r-1}
,w_{\epsilon r,-r-1}).
\end{align*}

Using this identity twice with $\phi=\phi^\sharp$ we obtain: 
\begin{equation}
a_r^4=(a_r\ sa_{-r})^2=\left(\frac{1}{ [r]!\
[T^{-1};r]_{(r-1)}}\right)^2.
\end{equation}
This gives the
second identity of \lemref{delicate} up to a sign. The correct sign is implied by the third
identity which we now prove.
\end{proof}

\subsection{} To verify the correct choice of sign for the third identity we
shall need some preliminary lemmas. Let $M^\prime$ denote the span of all
the weight subspaces of $M$ other than the highest weight space. Let $\d$
denote the projection of $M\otimes_R \mathcal F$ onto
$w_{0,-1}\otimes \mathcal F$ with kernel $M'\otimes_R \mathcal F$. Define constants
$c_\pm$ by the relations:
 $\delta(w_{r,r-1})\equiv c_+ w_{0,-1}\otimes u^{(n)}_k\ \mod \ M'\otimes
_R\mathcal F$  and
$\delta(w_{ -r,-r-1})\equiv c_- w_{0,-1}\otimes u^{(n)}_l\ \mod 
\ M'\otimes_R\mathcal F$ where $n-2p-1=r-1$ and $n-2l-1=-r-1$. For any integer
$t$ set $z_t=w_{0,-1}\otimes u^{(n)}_t$. In a similar fashion define the
projection
$\d^\vee$ of $M_F\otimes _R\mathcal F$ onto $w_{0,1}\otimes _R\mathcal F$ with
kernel 
$M^\vee\otimes _R\mathcal F$ and $M^\vee$ equal to the span of all weight
subspaces in $M_F$ for weights other than $m+1$.

\begin{lem}\label{lastlemma} For $r$ a non-negative integer
\begin{align}
\d^\vee(w_{r,r+1})&=c_+\ \frac{v^{(n-k)(k+1)}
[T^{-1};2k-n-1]_{(p)}}{ [T^{-1};p]_{(p)}} \ w_{0,1} \otimes u^{(n)}_p, \\
  -\frac{c_-}{ [l]!}&\equiv \frac{c_+}{ [p]!}\mod T-1, \notag \\
\text{and}
\quad
\delta(\bar\Psi w_{r,r+1})&\equiv c_-\ \frac{(-1)^{r+1}v^{2l(l-r+1)}}{ [r]!} z_l\mod
(T-1)\cdot P_r. \notag
\end{align}

\end{lem}
\begin{proof} Let $\mathcal F$ be an $\bu$-module that is $X$-admissible and finite dimensional, then we have the expansion (see
\cite[3.9.1]{MR96c:17022})
\begin{equation}
F^{-s}\otimes a \otimes e \ \mapsto \  \sum_{j \in  \mathbb N}(-1)^j
\qbinom{j+s-1}{ s-1}_i v^{-j(j+s)}\  F^{-j-s}a \otimes F^jK^{j+s}e \ ,
\end{equation}

Recall from \lemref{highestlowestwtvectors} 
\begin{equation}
w_{r,r-1}=c_+\sum_{0\le j
\le p}\frac{[n-p+j]_{(j)}v^{-j^2}}{T^{-j}
[T^{-1};j]_{(j)}} F^{(j)}w_{0,-1}\otimes u^{(n)}_{p-j}
\end{equation}
for $r=n-2p$.
 Then 
\begin{align*}
&F^{-1}w_{r,r-1} \\
&=c_+\sum_{j=0}^pv^{-j^2+(j+1)(n-2p+2j)}\frac{[n-p+j]_{(j)}}{T^{-j}
  [T^{-1};j]_{(j)}[j]![p-j]!} \\
&\hskip 100pt\times \Big(\sum_{l \in  \mathbb N}(-1)^l v^{-l(l+1)} [p-j+l]!  F^{j-l-1}w_{0,-1}\otimes u^{(n)}_{p-j+l}\Big).
\end{align*}
Set $u= w_{0,1}\otimes u^{(n)}_p$.
\begin{align*}
\d^\vee(F&^{-1}w_{r,r-1})
        =\frac{{c_+} v^{(n-p)(p+1)}}{ [T^{-1};p]_{(p)}}\
[T^{-1};2p-n-1]_{(p)}\
\cdot u .
\end{align*}  Here we have used the binomial identity given in \lemref{thirdbinomial}.  This implies 
\begin{align*}
w_{r,r+1}=F^{-1}w_{r,r-1}
        &=\frac{{c_+} v^{(n-p)(p+1)}}{ [T^{-1};p]_{(p)}}[T^{-1};2p-n-1]_{(p)}m_{r,r+1} .
\end{align*}  
Thus by \eqnref{needlater}
\begin{align*}
\bar\Psi(w_{r,r+1})
&=(-1)^{n-p} \frac{{sc_+} v^{-(n-p)(p+1)}}{[T;p]_{(p)}}[T;2p-n-1]_{(p)} \\
&\quad \quad \times \sum_{s=0}^{n-p}\frac{v^{s(1-r)+(n-p-s)(p+s+1)}[p+s]_{(s)}}{T^{-s}[T^{-1},s]_{(s)}}
   F^{(s)}w_{0,-1}\otimes u^{(n)}_{n-p-s}   \\
&=(-1)^{n-p} \frac{{sc_+} }{[T;p]_{(p)}}[T;2p-n-1]_{(p)} \\
&\quad \quad \times \sum_{s=0}^{n-p}\frac{v^{-s^2}[p+s]_{(s)}}{T^{-s}[T^{-1},s]_{(s)}}
   F^{(s)}w_{0,-1}\otimes u^{(n)}_{n-p-s}  \\
&\color{red}=(-1)^{n-p} \frac{{sc_+} }{[T;p]_{(p)}}[T;2p-n-1]_{(p)} w_{-r,-r-1}\color{black}
\end{align*}

In $P_r$ we know from the basis that $F^rw_{r,r-1}\equiv -w_{-r,-r-1}\mod
\ (T-1)\cdot P_r$.  Since the action of
$F$ commutes with $\d$,
we obtain: $[p]!\ {c_-}\equiv -[l]!\ {c_+} \mod\
T-1.$ 
For the last identity of the Lemma apply $\delta$ to the equation above. We get 
\begin{align*}
\delta(\bar\Psi w_{r,r+1})
&\equiv(-1)^{r+1}c_-\ \frac{1}{[r]!}\ z_l \end{align*}  
Here we use $r+p=l=n-p$ and $[T;r]\equiv [r]\mod (T-1)$.

\end{proof}

We now return to the proof of the congruence. Since $\Psi(w_{r,r-1})=-a_r
w_{-r,-r-1}$  we can calculate the constant $a_r$ as the ratio of
$\delta(\Psi(w_{r,r+1}))$ and $\delta( w_{-r,-r-1})$.  From \eqnref{lastlemma} we find the ratio
is congruent to $\frac{(-1)^{r+1}}{ [r]!}\mod T-1$. This completes the proof
of \lemref{delicate}.

\subsection{} Recall from (3.3) the category $_R\mathcal C_i$ and note
that any module $N$ in the category is the direct sum of generalized
eigenspaces for the Casimir element \cite{MR96d:17015} in the sense that 
$N=\sum N^{(\pm r)}$ where the sum is over $\mathbb N$ and $N^{(\pm r)}$
contains all highest  weight vectors in $N$ with weights $m+r-1$ and
$m-r-1$. Note that $N^{(\pm r)} $ need not be generated by its highest
weight vectors. The decomposition in \secref{decompositionsection}, $M\otimes_R\mathcal F\equiv \sum
P_i$ where the sum is over the nonnegative weights of $_R\mathcal F$ is such a
decomposition. In this case $(M\otimes _R\mathcal F)^{(\pm i)}= P_i$.  The
Casimir element $\Omega_0$ of $_R\mathbf U$ by 
$$
\Omega_0 =F E + \frac{vK-2+v^{-1} K^{-1}}{( v-v^{-1})^2}.
$$   
Let $N^{(r)}$ (resp. $N^{(-r)}$) denote the submodule of $N$ where the
Casimir element acts by the scaler 
$$ c(\l)=\frac{v^{r-1}T-2+v^{-r+1}T^{-1}}{( v-v^{-1})^2
}
\quad \text{\rm resp.}\quad  c(s\l)= \frac{v^{-r+1}T-2+v^{r-1} T^{-1}}{ (
v-v^{-1})^2}
$$

\subsection{} We now turn to the general case where $\mathcal F$ is a finite
dimensional $\bui$-module but not   necessarily irreducible. We extend the
definition of the
$s$-linear maps $\G$ and $\overline \G$ defined in \eqnref{ar} as follows. Decompose
$M\otimes_R \mathcal F$  into generalized eigenspaces for the Casimir
$(M\otimes_R \mathcal F)^{(\pm r)}$ and let $\mathfrak L^r$ denote the $m-r-1$ weight
subspace of
$(M\otimes_R \mathcal F)^{(\pm r)}$. Then set $\mathbb L=\sum \mathfrak L^r$. Decompose
$\mathcal F=\sum \mathcal F_{n_j}$  into irreducible $\rui$ modules.  Then $M\otimes_R
\mathcal F=\sum M\otimes_R \mathcal F_{n_j}$ and so we obtain $s$-linear  extensions
also denoted $\G$ and 
$\overline \G$ from $\mathbb L\cap (M\otimes_R \mathcal F_j)$ to all of $\mathfrak L$.  Set $(M\otimes _R\mathcal F)^{(+)}:
=\oplus_{r\in\mathbb N^*} M\otimes _R\mathcal F^{(r)}$

\begin{prop} Suppose $\phi$ is any invariant form on
$M\otimes_R \mathcal F$ with $\phi = \pm \phi^\sharp$. Let $\{w_j, j\in J\}$ be
an $R$-basis for the
 highest weight space of $(M\otimes_R \mathcal F)^{(0)}$ and
$\{u_i, i\in I\}$ a basis of weight vectors for the $E$-invariant weight
spaces of  weight $m+t$ for $t<-1$. Set
$M_j$ equal to the $\rui$-module generated by $w_j$ and $Q_l$ the
$\rui$-module generated by $(T-1)^{-1}(u_l+
\G(u_l))$. Then $M\otimes_R \mathcal F=\sum_jM_j\oplus \sum_lQ_l$ where each
$M_j\cong M(m)$ and  if $u_l$ has weight $m-t$, then $Q_l\cong P(m+t)$.
Moreover, if the basis vectors
$w_j$ and $u_l$ are $\phi$-orthogonal then the sum is an orthogonal sum of
$\rui$-modules.

\end{prop}
\begin{proof}
 Since   $M\otimes_R
\mathcal F=\sum M\otimes_R \mathcal F_{j}$ we may apply \corref{matrixcor} to each
summand to obtain an
$R$-basis of weight vectors $\{[T,0]\cdot x_i\}$ for $\mathfrak L\cap (M\otimes
_R\mathcal F)^{(+)}$ for which $\{[T,0]\cdot x_i,x_i+\ov\G x_i\}$ is a basis for
$\mathfrak L$ and the
$\ru$-module generated by $\mathfrak L$ is the direct sum of the submodules
generated by the vectors $x_i+\ov\G x_i$. If $x_i$ has weight
$m+t_i-1$ then $\ru\cdot(x_i+\ov\G x_i)\cong P(m+t_i)$ and
the intersection of this module with $\mathfrak L$ has $R$-basis $\{[T;0]\cdot
x_i,x_i+\ov\G x_i\}$. 

Now let $A$ denote the transition matrix from the basis 
$\{[T;0]\cdot u_i, i\in I\}$ to the basis $\{[T;0]\cdot x_i, i\in I\}$.
Then the determinant of $A$ is a unit of $R$. The block matrix
\begin{equation}
\begin{pmatrix}
 sA&0\\B&A
\end{pmatrix}
  \ ,\quad \text{for}\  [T;0]\cdot B=A-sA
\end{equation}
is the transition matrix from $\{u_i+\ov\G u_i ,[T;0]\cdot
u_i\}$ to $\{x_i+\ov\G x_i ,[T;0]\cdot
x_i\}$. The determinant is a unit and so the former set is a basis
of $\mathfrak L$. From this we conclude $M\otimes_R\mathcal F=\sum_jM_j\oplus
\sum_iQ_i$ and each $\mathfrak L\cap Q_i$ has $R$-basis $\{u_i+\ov\G u_i,[T;0]\cdot
u_i\}$. So if
the basis vectors
$w_j$ and $u_i$ are $\phi$-orthogonal then the sum $\mathfrak L=
\sum_j\mathfrak L\cap M_j\oplus
\sum_i\mathfrak L\cap Q_i$ is an
orthogonal sum. It follows that $M\otimes_R\mathcal F=\sum_jM_j\oplus
\sum_iQ_i$ is an orthogonal sum.
\end{proof}
%
%

\subsection{Diagonalizing forms.}

Fix a non-zero form $\phi$ on $A$, a free $R$-module of rank $n$. Choose an
integer $d_1$ with $\phi(A,A)=R\pi^{d_1}$ and choose vectors $a_1$ and
$a_1^\p$ with $\phi( a_1, a_1^\p)= \pi^{d_1}$

\begin{lem} Suppose the form $\phi$ is not zero. Set $A_1=\{a\in A |
\phi(a,a_1^\p)=0\}$ and
$A_1^\p=\{a\in A | \phi(a_1,a)=0\}$. Then we have direct sum
decompositions
$$R a_1\oplus A_1\cong A\cong R a_1^\p\oplus A_1^\p\ .$$
Moreover if $\phi$ is symmetric we may choose a unit of $R$, $u_1$ with
$u_1 a_1= a_1^\p$.
\end{lem}
\begin{proof} For $a\in A$ define
$$P(x)= x-\frac{\phi(x,a_1^\p)}{\phi(a_1,a_1^\p)} a_1\ ,$$
$$P^\p(x)= x-\frac{\phi(a_1,x)}{\phi(a_1,a_1^\p)} a_1^\p\ .$$
Then $P$ (resp. $P^\p$) is the projection of $A$ onto $A_1$ (resp.
$A_1^\p$). We have
$$x= P(x)+\frac{\phi(x,a_1^\p)}{\phi(a_1,a_1^\p)} a_1\ \quad,\quad
x=P^\p(x)+\frac{\phi(a_1,x)}{\phi(a_1,a_1^\p)} a_1^\p\ .$$
Since $a_1\notin A$ and $a_1^\p\notin A^\p$,
these are the
desired decompositions.
Now suppose that $\phi$ is symmetric. 
Consider the line between $a_1$
and $a_1^\p$,  $b_t=t a_1+
(1-t)a_1^\p$.
Then $\phi(b_t,b_t)=t^2 \phi(a_1,a_1)
+2t(1-t)\phi(a_1,a_1^\p)+(1-t)^2\phi(a_1^\p,a_1^\p)$. Since $t^2,
2t(1-t),$ and $(1-t)^2$ are linearly independent if we are not in
characteristic two  we can find an open set
of $t$ in the base field so that we have equality of ideals
$R\phi(b_t,b_t)=\phi(A,A)$.    Alternatively first suppose $\phi(a_1,a_1)=u_2\pi^{d_1}$ with $u_2$ a unit.
Then we have $\phi(a_1,u_2^{-1}a_1)=\pi^{d_1}$ and we have proved the remaining statement.  Similarly $\phi(a_1',a_1')=u_2'\pi^{d_1}$ with $u_2'$ a unit
leads to the same conclusion.  Lastly we suppose $\phi(a_1,a_1)=u_2\pi^{d_2}$ and $\phi(a_1',a_1')=u_2'\pi^{d_2'}$ with $d_2,d_2'>d_1$.
Then
\begin{align*}
\phi(b_t,b_t)&=t^2 \phi(a_1,a_1)
+2t(1-t)\phi(a_1,a_1^\p)+(1-t)^2\phi(a_1^\p,a_1^\p) \\
&=(2t(1-t)+t^2 u_2\pi^{d_2-d_1}+(1-t)^2u_2'\pi^{d_2'-d_1})\pi^{d_1} \\
&=u_t\pi^{d_1} 
\end{align*}
where $u_t$ is a unit for $t\neq 0,1$.

 From this we find a unit $u_1\in R$ such that
$ \phi(b_t,u_1 b_t) =u_1 \phi(b_t,b_t)=\pi^{d_1}$.
\end{proof}

Directly as a corollary to this lemma we have

\begin{cor}  There exist integers $m\leq n$ and $d_i, 1\leq
i\leq m$ and two bases for $A$, $\{a_i|1\leq i\leq n\}$ and $\{a^\p_i|1\leq i\leq n\}$ with the
following property
$$
\phi(a_i,a_j^\p)=\begin{cases}\delta_{i,j}\ \pi^{d_i} & \ i\leq m,  \\
0& \ m< i\leq n.\end{cases}$$

Moreover if $\phi$ is symmetric we may choose the bases so that
$a_i^\p=u_i a_i$ for some choice of units $u_i, 1\le i\le n$.
\end{cor}

\begin{proof}  If the form is zero then any two bases will suffice. So
assume $\phi$ is not zero and apply the lemma obtaining vectors $a_1$
, $a_1^\p$ and submodules  $A_1$ , $A_1^\p$. Proceed inductively and
assume for some
$t$ with $ 1\le t\le n-1$ and $m'\leq m$, we have chosen integers $d_i$ and vectors
$\{a_i\}$ and
$\{a_i^\p\}$ which satisfy the conditions 
$$\phi(a_i,a_j^\p)=\delta_{i,j}\ \pi^{d_i} \ \ i\leq t$$
$$ \phi(a_i,a_j^\p)=\delta_{i,j}\ \pi^{d_i} \ \ i\leq m' $$
$$\phi(a_i,a_j^\p)=0\  \ ,\ m< i\leq t \ .$$
 Also assume that for
  $A_t=\{a\in A |
\phi(a,a_i^\p)=0, 1\le i\le t\}$ and
$A_t^\p=\{a\in A | \phi(a_i,a)=0 , 1\le i\le t\}$, we have direct
sum decompositions
$$\sum_{1\le i\le t} R a_i\oplus A_t\cong A\cong \sum_{1\le i\le t} R
a_i^\p\oplus A_t^\p\ .$$
Now choose $a_{t+1}$ and $a_{t+1}^\p$ so that
$\phi(a_{t+1},a_{t+1}^\p)=\phi(A_t,A_t^\p)$. Then applying the lemma
gives the following decomposition: $R a_{t+1}\oplus A_{t+1}\cong A_t$, 
$ R a_{t+1}^\p\oplus A_{t+1}^\p\ \cong A_t^\p$. This completes the
inductive step which in turn proves the corollary.

\end{proof}

\begin{cor}  Suppose $\phi$ is symmetric. There exist integers
$m\leq n$ and
$d_i$ and units $u_i$ of $R, 1\leq i\leq m$ and a basis for $A$, $\{a_i\}$ ,
such that
$\phi$ is represented by the diagonal $n\times n$ matrix $S$ with entries:
$$S_{i,i}=u_i\ \pi^{d_i}\quad  \text{for}\quad i\leq m,\quad \text{and }\quad S_{i,i}=0\quad i>m  $$
$$ S=\text{diag}(u_1\pi^{d_1},\dots ,u_m\pi^{d_m},0,\dots, 0) $$
\end{cor}

Note that we cannot get rid of the units in the symmetric form case since doing
so would require a square root of each unit $u_i$. These square roots may not
lie in the ring $R$.

\begin{prop}
Suppose $\phi$ is any invariant  symmetric  form on
$M\otimes_R \mathcal F$ with $\phi = \pm \phi^\sharp$. Then $M\otimes_R \mathcal F$
admits an orthogonal decomposition with each summand an indecomposable
$\ru$-module and isomorphic to $M$ or some $P(m+t)$ for $t\in\mathbb N^*$.
\end{prop}

\begin{proof}   Since $R$ is a discrete valuation ring we may choose
an orthogonal $R$-basis for the free $R$-module $\mathfrak L\cap (M\otimes_R\mathcal F)^{(+)}$.
\end{proof}

\section{Filtrations}

\subsection{} We continue with the notation of the previous section. So
$\phi$ is an invariant form on
$M\otimes_R \mathcal F=\sum_iP_i$. For any $R$-module $B$ set $\overline
B=B/(T-1)\cdot B$ and for any  filtration $B=B_0\supset B_1\supset
...\supset B_r$, let $\overline B=\overline B_0\supset \overline B_1\supset
...\supset
\overline B_r$ be the induced filtration of $\overline B$, with $\overline
B_i=( B_i+(T-1)\cdot B)/ (T-1)\cdot B$.
Now $\phi$ induces a filtration on $M\otimes_R \mathcal F$ by
\begin{equation}\label{filtration}
(M\otimes_R \mathcal F)^i=\{v\in  M\otimes_R \mathcal F| \phi(v,M\otimes_R \mathcal
F)\subset (T-1)^i\cdot R\}.
\end{equation}

\subsection{}
 Let notation be as in (8.3) with $P=P_r$ and $\mathfrak L$ equal to the
$m-r-1$ weight subspace of $P$. Suppose $P=A_0\supset A_1\supset ...\supset
A_t=0$ is a filtration. Then since the $U$-module $P/(T-1)\cdot P$ contains only the
two proper subspaces $M(\pm r)$ we can
choose constants $a\ge b \ge c$,  so that
\begin{align}
\overline P = \overline A_0&=\cdots=\overline A_a \\
&\cup  \notag\\
M(r)\cong \overline A_{a+1}&=\cdots =\overline A_b \notag  \\
 &\cup \notag \\
 M(-r)\cong \overline A_{b+1}&=\cdots = \overline A_c \notag \\
&\cup   \notag\\
\overline A_{c+1}&=0\ .\notag
\end{align}
In this case we say that the filtration is of {\it type} $(a,b,c)$.   When $a=b=c$ we use the convention $\overline P=M(0)$. 

Recall from \secref{decompositionsection} the bases $w_{\pm r,-r-1}$ of the $m-r-1$ weight space of $M(\pm r)$
and the basis $\{ z_{r,-r-1}= w_{r,-r-1}+w_{-r,-r-1}\ ,\ [T;0] w_{r,-r-1} \}$ of the
$m-r-1$ weight space of $P$. Fix $r$ and for convenience set
\begin{equation}
w_+=w_{r,-r-1}\ ,\ w_-=w_{-r,-r-1} \ ,\ z=w_+ +w_- \ .
\end{equation}

Let $\phi_\pm$ denote the Shapovalov form on the Verma module $M_{\pm r}$ and
normalized by the identities
\begin{equation}
\phi_+(w_+,w_+)=1  \quad \quad  \phi_-(w_-,w_-)=1\ .
\end{equation}
Now choose constants $b_\pm\in R$ with $\phi|_P=b_+\phi_+ +b_- \phi_-\ .$ Then the
restriction of $\phi$ to $\mathfrak L$ is given by the matrix $\mathcal M$ with respect
to the basis $\{z,[T;0]w_+\}$
$$\mathcal M=
\begin{pmatrix}
\phi(z,z)& \phi(z, [T;0] w_+)\\ \phi([T;0] w_+,z) & \phi([T;0] w_+,[T;0] w_+)
\end{pmatrix}
=
 \begin{pmatrix}
b_+ +b_- & \  [T;0] b_+ \\  \ [T;0] b_+ & \ [T;0]^2 b_+
\end{pmatrix}
$$

Define the order of elements $a\in  \mathcal K$ by $ord(a)=n$ if $a\in R(T-1)^n$ and
$a\notin R(T-1)^{n+1}$. For a matrix define the order to be the minimum of the orders
of the matrix entries. For any matrices $\mathcal A$ and $\mathcal B$ with entries in
$R$, $ord (\mathcal A \mathcal B) \ge ord \mathcal A$ and if $det \mathcal A$ is a
unit then $ord (\mathcal A \mathcal B) = ord (\mathcal B)$.

The type of the filtration on $P$ can now be determined easily by the constants
$b_\pm$. The result separates into three cases.

\begin{lem}\label{formonP} \begin{enumerate}[(a).]
\item Suppose $ord(b_+)=ord(b_-)< ord(b_++b_-)$. Then the
filtration on $P$ is of type $(a,a,a)$ with $a= ord(b_+)+1$.
\item Suppose $ord(b_++b_-)=ord(b_+)\leq  ord(b_-)$. Then the
filtration on $P$ is of type
$( ord(b_+), ord(b_+)+1, ord(b_-)+2 )$.
\item Suppose $ord(b_++b_-)=ord(b_-)\leq  ord(b_+)$. Then the
filtration on $P$ is of type $( ord(b_-), ord(b_+)+1, ord(b_+)+2 )$.
\end{enumerate}

\end{lem}
\begin{proof} Let $(a,b,c)$ designate the type of filtration on $P$ induced by the
form $\phi$.  Diagonalize $\mathcal M$ as follows. Choose invertible $R$-valued
$2\times 2$ matrices
$\mathcal U$ and $\mathcal V$ and integers $d_1\leq d_2$ so that
$$\mathcal U \mathcal M \mathcal V = \begin{pmatrix} (T-1)^{d_1} & \ \\  \ &
(T-1)^{d_2}
\end{pmatrix} = \mathcal D\ .  $$
Then $ord(\mathcal M)= ord(\mathcal D)=d_1$ and $d_2=ord(det(\mathcal D))-d_1=
ord(det(\mathcal M))-ord(\mathcal M)$. We conclude $a=d_1= ord(\mathcal M)$ and so
from the form of $ \mathcal M$ above we get the formulas for $a$ in the three cases.
Similarly
$c=d_2= ord(det(\mathcal M))-ord(\mathcal M)$ and in the three cases of the lemma this
translates to the formulas in the lemma.

Next we determine $b$. The highest weight space of $P$ is a free rank one $R$-module.
So $b= ord(\phi([T;0] w_{r,r-1} , [T;0] w_{r,r-1}))$. Then for some unit $u\in R$,
$$\phi([T;0] w_{r,r-1}\ ,\ [T;0] w_{r,r-1})=\phi(u E^r z\ ,\ [T;0] w_{r,r-1})=
\phi(u z\ ,\ F^r [T;0] w_{r,r-1}) = \ $$
$$\phi(u z\ ,\  [T;0] w_{r,-r-1}) = u [T;0] b_+\ .$$
So the form restricted to the highest weight space has order $ord(b_+)+1$. This proves
$b= ord(b_+)+1$ and completes the proof of the lemma.

\end{proof}

\section{Examples}    In this section we will roughly follow the notation in \cite{MR94m:17016}.  
A pair $(\Pi,(\enspace,\enspace))$ where $\Pi$ is a finite set and $(\enspace,\enspace)$ denotes a
symmetric bilinear form on the free abelian group $\mathbb Z[\Pi]$ with
values in $\mathbb Z$ is called a {\it Cartan datum} if  
\begin{equation}
(\alpha, \alpha)\in\{2,4,6,\dots\}\text{ for  any } \alpha \in \Pi ;\quad 2
\frac{(\alpha,\beta)}{ (\alpha,\alpha)}\in\{0,-1,-2,\dots\}  \text{ for any }  \alpha\neq \beta
\in \Pi.,
\end{equation}
If $P\in \Bbb Q(v)$ is a
rational function then $P_\alpha$ denotes $P(v_\alpha)$ where $v_\alpha=v^{(\alpha,\alpha)/2}$\label{Gaussintegers}.
\subsection{Verma modules for the quantum group $U_v(\mathfrak{sl}(3))$.} Let 
\[
\Phi=\{\pm \alpha,\pm\beta,\pm(\alpha+\beta)\}
\]
 be the root system of $\mathfrak{sl}_3(\mathbb C)$ with $\Pi=\{\alpha,\beta\}$ a set of simple roots, $W$ the Weyl group, and $(\enspace,\enspace)$  the unique  $W$-invariant form defined on $\Phi$ with $(\gamma,\gamma)=2$ for all $\gamma\in \Phi$.  Let $\varpi_\gamma$ denote a fundamental weight with respect to $\Pi$, $\gamma\in \Pi$, and let the weight lattice of $\Phi$ be $\Lambda=\sum_{\gamma\in\Pi}\mathbb Z\varpi_\gamma$\label{Lambda}. 
 
  \color{black} The quantum enveloping algebra $U_v(\mathfrak{sl}_3(\mathbb C))$ is
 defined to be the associative algebra over $\mathbb Q[v,v^{-1}]$ with generators $E_\gamma$, $F_\gamma$, $K_\gamma$ and $K^{-1}_\gamma$, ($\gamma\in \Pi$) subject to the relations 
\begin{gather*}
K_\gamma K^{-1}_\gamma=K_\gamma^{-1}K_\gamma=1, \quad K_\gamma K_\nu =K_\nu K_\gamma \tag{R1}\\
K_\gamma E_\nu K^{-1}_\gamma=v^{(\gamma,\nu)}E_\nu ,\quad
K_\gamma F_\nu K^{-1}_\gamma=v^{-(\gamma,\nu)}F_\nu \tag{R2} \\
[E_\gamma,F_\nu]=\delta_{\gamma,\nu} \frac{K_\gamma-K_{-\gamma}}{v_\gamma-v^{-1}_\gamma} \tag{R3}\\
E_\gamma^2 E_\nu -[2] E_\gamma E_\nu E_\gamma+ E_\nu E_\gamma^2=0  \quad \text{ for }\gamma\neq \nu\tag{R4} \\
F_\gamma^2 F_\nu -[2] F_\gamma F_\nu F_\gamma+ F_\nu F_\gamma^2=0 \quad \text{ for }\gamma\neq \nu   
\end{gather*}
where $\gamma,\nu\in \Pi$.
One also sets $\rho=\varpi_\alpha+\varpi_\beta=\alpha+\beta$, $U=U_v(\mathfrak{sl}(3)))$ and let $U^+$ (resp. $U^-$) denote the subalgebra of $U$ generated by $E_\nu$ (resp. $F_\nu$) with $\nu\in \Phi$.   Moreover let $U_v(\mathfrak a)$\label{ua} denote the subalgebra generated by $E_\alpha,F_\alpha, K_\gamma$ with $\gamma\in \Pi$.  We will now consider $M(-\varpi_\alpha)$ which has two Verma submodules $M(-2\varpi_\beta)$ and $M(-2\varpi_\alpha-\varpi_\beta)$ due to the fact that $s_\beta(-\varpi_\alpha+\rho)-\rho=-\varpi_\alpha-\beta=-2\varpi_\beta$ and $s_{\alpha}(-2\varpi_\beta+\rho)-\rho=-2\varpi_\alpha-\varpi_\beta$.

  The subalgebra $U^-$ generated by $F_\gamma$, $\gamma\in\Pi$, has a basis of the form 
 $$
T_\alpha T_\beta (F_\alpha^{(m)})T_\alpha(F^{(n)}_\beta)F_\alpha^{(p)} =T_\alpha \left(T_\beta (F_\alpha^{(m)})F^{(n)}_\beta\right)F_\alpha^{(p)} 
$$
where $m,n,p\in\mathbb N$ and $T_\alpha=T''_{1,1}, T_\beta=T''_{2,1}$ (in Lusztig's notation, see \cite{MR94m:17016} and \cite[Theorem 8.24]{MR96m:17029}).
By \cite[8.16.(6)]{MR96m:17029} $T_\alpha T_\beta(F_\alpha)=F_\beta$, so we get
\begin{equation*}
T_\alpha T_\beta (F_\alpha^{(m)}K_\alpha^m)= F_\beta^{(m)}K_\beta^{m},
\end{equation*}
and
\begin{equation}\label{calc1}
T_\beta(F_\alpha K_\alpha)=\left(F_\alpha F_\beta -vF_\beta F_\alpha\right)K_{\beta+\alpha}.
\end{equation}

On the other hand
\begin{align*}
\text{ad}\,(F_\alpha)(F_\beta K_\beta) =(F_\alpha F_\beta -vF_\beta F_\alpha)K_{\alpha+\beta}
\end{align*}


%

 In the notation of section \secref{QCG decomposition}, we have under the adjoint action on $U$,
\begin{equation}\label{F1}
 u^{(1)}=F_\beta K_\beta ,\quad u^{(1)}_1=\text{ad}\,(F_\alpha)(F_\beta K_\beta) 
\end{equation}
 and $\text{ad}\,(F_\alpha^{(2)})(F_\beta K_\beta)=0$ by \cite[Lemma 4.18]{MR96m:17029} so that the set $\{u^{(1)}$, $u^{(1)}_1\}$ spans a copy of $\mathcal F_1$ inside $U^-$.

Let $M(\lambda)$ denote the Verma module for $U_\nu(\mathfrak{sl}(3))$ with highest weight $\lambda\in \Lambda$  defined by 
$$
M(\lambda):=U/\left( \sum_{\nu\in \Pi}UE_\nu +\sum_{\nu \in\Pi} U(K_\nu- T^{\nu\cdot\nu/2}v^{(\lambda, \nu)})\right).
$$
Let $\mathbf 1$ denote the image of $1$ in the quotient $M(\lambda)$.  Set
$$
M(\lambda)^\nu:=\{m\in M(\lambda)\,|\, K_\gamma m=T^{\nu\cdot\nu/2} v^{(\gamma,\lambda)}\text{ for all }\gamma\in\Pi\}.
$$
The following is a quantum analogue of a result Shapovalov's: 
\begin{prop}[{\cite[Proposition 19.1.2]{MR94m:17016}}]  For any $\lambda\in \Lambda$ there exists a unique symmetric bilinear form $M(\lambda)\times M(\lambda)\to R$ such that 
\begin{enumerate}
\item $(\mathbf 1,\mathbf 1)=1$;
\item $(ux,y)=(x,\varrho(u)y)$ for all $x,y\in M(\lambda)$ and $u\in\bu$.
\end{enumerate}
Moreover $(x,y)=0$ if $x\in M(\lambda)^\nu$ and $y\in M(\lambda)^{\nu'}$ with $\nu\neq \nu'\in \Lambda$ .
\end{prop}

Since the Shapovalov form is $\varrho$-invariant it must be an induced form when restricted to a $U_\nu(\mathfrak a)$-summand isomorphic to $ U_\nu(\mathfrak a)\mathbf 1\otimes\mathcal F_n$.  Here $\mathcal F_n$ is any ad-invariant irreducible summand of $\bu$.  To illustrate what is going on in the paper we use the copy of $\mathcal F_1$ given in \eqnref{F1}, taking into account that we need to use the $_f\mathcal R^{-1}:\mathcal F_n \otimes U_\nu(\mathfrak a)\mathbf 1\to U_\nu(\mathfrak a)\mathbf 1\otimes\mathcal F_n$ (we will use $f$ as defined in \eqnref{fexample}).  More precisely let us determine $\boldsymbol{\beta}:\mathcal F_1\otimes \mathcal F_1^{\rho_1}\to U_v(\mathfrak a)$ satisfying 
\begin{equation}\label{shap}
(um,u'n)=\chi_{\boldsymbol\beta,\phi_M}({_f\mathcal R}^{-1}(u\otimes m), {_f\mathcal R}^{-1}(u'\otimes n))
\end{equation}
where $u,u'\in \mathcal F_1$, $m,n\in M: =U_\nu(\mathfrak a)\mathbf 1$ and $\phi_M$ is the Shapovalov form on the Verma $U_v(\mathfrak a)$-module of highest weight $Tv^{-1}$ with respect to $K_\alpha$ ($T$ with respect to $K_\beta$) and is normalized so that $\phi_M(\mathbf 1,\mathbf 1)=1$.  Note that the linear map $\boldsymbol{\beta}$ has a bold font to distinguish it from the root $\beta$.
Since
$$
\mathcal F_1\otimes \mathcal F_1^{\rho_1}\cong\mathcal F_{0}\oplus \mathcal F_2,
$$
the problem then is to find $r_i\in\mathbb C[T]$ with
$$
\boldsymbol{\beta}=r_0\boldsymbol{\beta}^{1,1}_0+r_2\boldsymbol{\beta}^{1,1}_2
$$
where 
$$
\boldsymbol{\beta}^{1,1}_{2r}(u^{(a)})=\delta_{2r,a}E_\alpha^{(r)}K_\alpha^{-r}.
$$
First of all (for $\mathbf 1$ the highest weight vector) using \eqnref{antiauto}, 
$$
(F_\beta K_\beta \mathbf 1,F_\beta K_\beta \mathbf 1)=T^2(\mathbf 1,\brho(F_\beta)F_\beta\mathbf 1)=T^2v (\mathbf 1, K_\beta^{-1} E_\beta F_\beta \mathbf 1)=vT[T;0].
$$
Now we use $f$ as \eqnref{fexample}
\begin{align*}
_f\mathcal R^{-1}(F_\beta K_\beta \otimes \mathbf 1)
&=f(1,-1)^{-1}\mathbf 1 \otimes F_\beta K_\beta.
\end{align*}
which has as a consequence 
\begin{align*}
(F_\beta K_\beta\mathbf 1 ,F_\beta K_\beta\mathbf 1 )&=\chi_{\boldsymbol\beta,\phi_M}(\mathbf s\circ \Pi_f^{-1}(F_\beta K_\beta\otimes \mathbf 1 ),\mathbf s\circ \Pi_f^{-1}(F_\beta K_\beta\otimes \mathbf 1)) \\
&=f(1,-1)^{-2}\phi_M(\mathbf 1, \rho_1\boldsymbol\beta(F_\beta K_\beta\otimes  F_\beta K_\beta)\mathbf 1).
\end{align*}

So by \corref{firstbetacor} we get
\begin{align}\label{relation1}
f(1,-1)^{-2}\rho_1\boldsymbol{\beta}(F_\beta K_\beta\otimes F_\beta K_\beta)\mathbf 1
&=f(1,-1)^{-2}\left(\frac{r_0}{[2]}+\frac{v r_2}{[2]}[T,-1]\right)\mathbf 1 
=vT[T;0]\mathbf 1 
\end{align}

Next by \eqnref{calc1}
\begin{align*}
(\text{ad}F_\alpha(F_\beta K_\beta)\mathbf 1,\text{ad}F_\alpha(F_\beta K_\beta)\mathbf 1)
&=vT [T;0]+v^2T^3[T;-1].
\end{align*}
On the other hand \corref{firstbetacor} gives us
\begin{align}
f(-1,-1)^{-2}\rho_1\left(r_0\boldsymbol{\beta}^{1,1}_0+r_2\boldsymbol{\beta}^{1,1}_2\right)&(\text{ad}F_\alpha(F_\beta K_\beta)\otimes \text{ad}F_\alpha(F_\beta K_\beta))\mathbf 1  
\label{relation2} \\
&= 
\left(vT [T;0]+v^2T^3[T;-1]\right) \mathbf 1\notag
\end{align}

Using the equations \eqnref{relation1} and \eqnref{relation2} we finish the determination of $\boldsymbol\beta$:
$$
 \boldsymbol{\beta}=f(1,-1)^{2}\left(T^2(T+T^{-1})[T;0]+vT^5[T;-1])\boldsymbol{\beta}^{1,1}_0-T^2(\{1\}  [T;0] +T^3)\boldsymbol{\beta}^{1,1}_2\right).
$$

Here however we have only determined one $\boldsymbol \beta$ for one $U_\nu(\mathfrak a)$-summand of $M(-\varpi_\alpha)$.  In future we plan to investigate the other $\boldsymbol\beta$ that appear.

\subsection{Verma Modules for $U_\nu(\mathfrak{sp}(4))$}
Let 
\[
\Phi=\{\pm \alpha,\pm\beta,\pm(\alpha+\beta),\pm (2\alpha+\beta)\}
\]
 be the root system of type $B_2$ with $\Pi=\{\alpha,\beta\}$ a set of simple roots, $\alpha$ short, $
 \beta$ long, $W$ the Weyl group, and $(\enspace,\enspace)$  the unique  $W$-invariant form defined on $\Phi$ with $(\alpha, \beta)=-2$, $(\alpha,\alpha)=2$ and $(\beta,\beta)=4$,
 so that
 $
 \langle \alpha,\check\beta\rangle
 =-1,
  \langle \beta,\check\alpha\rangle
  =-2.
 $
 Let $\varpi_\gamma$ denote a fundamental weight with respect to $\Pi$, $\gamma\in \Pi$, and let the weight lattice of $\Phi$ be $\Lambda=\sum_{\gamma\in\Pi}\mathbb Z\varpi_\gamma$.   Then $q_\alpha=q$ and $q_\beta=q^2$.  
 
 The quantum enveloping algebra $U_v(\mathfrak{sp}(4))$ is
 defined to be the associative algebra over $\mathbb Q[v,v^{-1}]$ with generators $E_\gamma$, $F_\gamma$, $K_\gamma$ and $K^{-1}_\gamma$, ($\gamma\in \Pi$) subject to the relations (R1)-(R3) but with a different Serre relation:
\begin{gather*}
E_\gamma^3 E_\nu -[3]_\gamma E_\gamma^2 E_\nu E_\gamma+[3]_\gamma E_\gamma E_\nu E_\gamma^2-E_\nu E_\gamma^3=0  \quad \text{ for }\gamma\neq \nu \\
F_\gamma^3 F_\nu -[3]_\gamma F_\gamma^2 F_\nu F_\gamma+[3]_\gamma F_\gamma F_\nu F_\gamma^2-F_\nu F_\gamma^3=0   \quad \text{ for }\gamma\neq \nu 
\end{gather*}
where $\gamma,\nu\in \Pi$.
$U=U_v(\mathfrak{sp}(4)))$ and let $U^+$ (resp. $U^-$) denote the subalgebra of $U$ generated by $E_\nu$ (resp. $F_\nu$) with $\nu\in \Phi$.   Moreover let $U_v(\mathfrak a)$ denote the subalgebra generated by $E_\alpha,F_\alpha, K_\gamma$ with $\gamma\in \Pi$.  Let $\rho_1:\bold U\to \bold U$ be the algebra isomorphism determined by
the assignment 
\begin{equation}
\rho_1(E_\gamma)=-v_\gamma F_\gamma,\quad \rho_1(F_\gamma)=-v_\gamma^{-1}E_\gamma,\quad
\rho_1(K_\gamma)=K_{\gamma}^{-1}
\end{equation}
for all $\gamma\in\Pi$.   Define also an algebra
anti-automorphism $\brho:\bold U\to \bold U$ by 
\begin{equation} \label{antiauto2}
\brho(E_\gamma)=v_\gamma K_\gamma F_\gamma,\quad \brho(F_\gamma)=v_\gamma 
K_{\gamma}^{-1}E_\gamma,\quad \brho(K_\gamma)=K_\gamma. 
\end{equation}
These maps are related through the antipode $S$ of  $\bold U$
by $\brho=\rho_1S$ .

We will now consider (a $U_\nu(\mathfrak a)$-submodule of) $M(-\varpi_\alpha+\varpi_\beta)$ which has three Verma submodules $M(-\varpi_\alpha+\varpi_\beta-2\beta)$, $M( -\varpi_\alpha+\varpi_\beta-4\alpha-2\beta)$, and $M(-\varpi_\alpha+\varpi_\beta-4\alpha-4\beta)$ due to the fact that 
$$
s_\beta(-\varpi_\alpha+\varpi_\beta+\rho)-\rho
=-\varpi_\alpha+\varpi_\beta-2\beta,
$$
$$
s_{\alpha}(-\varpi_\alpha+\varpi_\beta-2\beta+\rho)-\rho
=-\varpi_\alpha+\varpi_\beta-4\alpha-2\beta,
$$
$$
s_{\beta}(-\varpi_\alpha+\varpi_\beta-2\beta-4\alpha+\rho)-\rho
=-\varpi_\alpha+\varpi_\beta-4\alpha-4\beta.
$$


Note that 
\begin{align}
\text{ad}\,(F_\alpha)(F_\beta K_\beta) &=(F_\alpha F_\beta -v^2F_\beta F_\alpha)K_{\alpha+\beta}\text{ and }\label{calc4}  \\
\text{ad}\,(F_\alpha^{(2)})(F_\beta K_\beta) 
&=\Big(F_\alpha^{(2)} F_\beta -vF_\alpha F_\beta F_\alpha +v^2F_\beta F_\alpha^{(2)}\Big)K_{2\alpha+\beta}.
\end{align}

\color{black}
 In the notation of section \secref{QCG decomposition}, we have under the adjoint action on $U$,
\begin{equation}\label{F2}
 u^{(2)}=F_\beta K_\beta ,\quad u^{(2)}_1=\text{ad}\,(F_\alpha)(F_\beta K_\beta) ,\quad u^{(2)}_2=\text{ad}\,(F_\alpha^{(2)})(F_\beta K_\beta) 
\end{equation}
 and $\text{ad}\,(F_\alpha^{(3)})(F_\beta K_\beta)=0$ by \cite[Lemma 4.18]{MR96m:17029} so that the set $\{u^{(2)}, u^{(2)}_1, u^{(2)}_2\}$ spans a copy of $\mathcal F_2$ inside $U^-$.
Let us determine $\boldsymbol{\beta}:\mathcal F_2\otimes \mathcal F_2^{\rho_1}\to U_v(\mathfrak a)$ satisfying 
\begin{equation}
(um,u'n)=\chi_{\boldsymbol\beta,\phi_M}({_f\mathcal R}^{-1}(u\otimes m), {_f\mathcal R}^{-1}(u'\otimes n))
\end{equation}
where $u,u'\in \mathcal F_1$, $m,n\in M: =U_\nu(\mathfrak a)\mathbf 1$ and $\phi_M$ is the Shapovalov form on the Verma $U_v(\mathfrak a)$-module of highest weight $v^{-1}T$ with respect to $K_\alpha$ ($v^2T$ with respect to $K_\beta$  - this is because $(\varpi_\beta,\beta)=2$ and $(\varpi_\alpha,\alpha)=1$) and is normalized so that $\phi_M(\mathbf 1,\mathbf 1)=1$.  
Since
$$
\mathcal F_2\otimes \mathcal F_2^{\rho_1}\cong\mathcal F_{0}\oplus \mathcal F_2\oplus
\mathcal F_4
$$
the problem then is to find $r_i\in\mathbb C[T]$ with
$$
\boldsymbol{\beta}=r_0\boldsymbol{\beta}^{2,2}_0+r_2\boldsymbol{\beta}^{2,2}_2
+r_4\boldsymbol{\beta}^{2,2}_4
$$
where 
$$
\boldsymbol{\beta}^{2,2}_{2r}(u^{(a)})=\delta_{2r,a}E_\alpha^{(r)}K_\alpha^{-r}.
$$
First of all (for $\mathbf 1$ the highest weight vector) using \eqnref{antiauto2}, 
$$
(F_\beta K_\beta \mathbf 1,F_\beta K_\beta \mathbf 1)=v^2T^4(\mathbf 1,\brho(F_\beta)F_\beta\mathbf 1)=T^4v^4 (\mathbf 1, K_\beta^{-1} E_\beta F_\beta \mathbf 1)=v^3T^2[T^2;1]_\beta.
$$
Next
from \corref{firstbetacor} (with $m=n=2$) we get
\begin{align*}
\Big(r_0\boldsymbol{\beta}^{1,1}_0&+r_2\boldsymbol{\beta}^{2,2}_2+r_4\boldsymbol{\beta}^{2,2}_4\Big)((\text{ad}\,F_\alpha^{(i)})F_\beta K_\beta\otimes (\text{ad}\,F_\alpha^{(i)})F_\beta K_\beta) \\
&=(-1)^iv^{2i(i-3)}{\qbinom 2 i}^2  \left(\sum_{k=0}^2r_{2k}
   \frac{\qbinom{T;-1}{k}}{||u^{(2k)}||^2\left[\begin{matrix} 2k
  \\  k \end{matrix}\right]} \left[\begin{matrix} 2 & 2 & 2k
  \\ i & {2-i} & k \end{matrix}\right]
\right)
\end{align*}

Now we use $f$ as \eqnref{fexample}
\begin{align*}
_f\mathcal R^{-1}(F_\beta K_\beta \otimes \mathbf 1)
&=f(2,-1)^{-1}\mathbf 1 \otimes F_\beta K_\beta.
\end{align*}
which has as a consequence 
\begin{align*}
(F_\beta K_\beta\mathbf 1 ,F_\beta K_\beta\mathbf 1 )&=\chi_{\boldsymbol\beta,\phi_M}(\mathbf s\circ \Pi_f^{-1}(F_\beta K_\beta\otimes \mathbf 1 ),\mathbf s\circ \Pi_f^{-1}(F_\beta K_\beta\otimes \mathbf 1)) \\
&=f(2,-1)^{-2}\phi_M(\mathbf 1, \rho_1\boldsymbol\beta(F_\beta K_\beta\otimes  F_\beta K_\beta)\mathbf 1).
\end{align*}

From this simple calculation and \corref{firstbetacor} we have
\begin{align}\label{relation3}
f(2,-1)^{-2}\rho_1\boldsymbol{\beta}(F_\beta K_\beta\otimes F_\beta K_\beta)\mathbf 1
&=f(2,-1)^{-2}\left(\frac{r_0}{[3]}+
\frac{ v^2\,r_2}{[4]}[T;-1]+\frac{v^4 [2]\,r_4}{[4][3]}\qbinom{T;-1}{2}\right)\mathbf 1 \\
&=v^4T^2[T^2;1]_\beta\mathbf 1 
\end{align}

Next by \eqnref{calc4}
\begin{align*}
(\text{ad}F_\alpha(F_\beta K_\beta)\mathbf 1,&\text{ad}F_\alpha(F_\beta K_\beta)\mathbf 1)
=((F_\alpha F_\beta -v^2F_\beta F_\alpha )K_{\alpha+\beta}\mathbf 1,(F_\alpha F_\beta -v^2F_\beta F_\alpha )K_{\alpha+\beta}\mathbf 1)\\
&=vT^2[2][T^2;1]_\beta +v^6T^5[T;-1]. 
\end{align*}
On the other hand
\begin{align*}
(\text{ad}F_\alpha(F_\beta K_\beta)\mathbf 1 ,\text{ad}F_\alpha(F_\beta K_\beta)\mathbf 1 )
&=f(0,-1)^{-2}\phi_M(\mathbf 1, \rho_1\boldsymbol\beta(\text{ad}F_\alpha(F_\beta K_\beta)\otimes  \text{ad}F_\alpha(F_\beta K_\beta))\mathbf 1).
\end{align*}
Again from \corref{firstbetacor} we obtain 

\begin{align}\label{relation4}
f(2,&-1)^{-2}v^2T^{-2}\rho_1\left(r_0\boldsymbol{\beta}^{2,2}_0+r_2\boldsymbol{\beta}^{2,2}_2+r_4\boldsymbol{\beta}^{2,2}_4\right)(\text{ad}F_\alpha(F_\beta K_\beta)\otimes \text{ad}F_\alpha(F_\beta K_\beta))\mathbf 1  
\\
&=f(2,-1)^{-2}v^2T^{-2}\left(\frac{r_0v^{-2}[2]}{[3]}+\frac{r_2v^{-1} 
   [2]\{1\}}
   {\left[4\right]} \qbinom{T;-1}{1}-r_4v^{-1}[2]\qbinom{4}{2}^{-1}\qbinom{T;-1}{2}\right)\mathbf 1
\notag \\ 
&= 
\left(vT^2[2][T^2;1]_\beta +v^6T^5[T;-1]\right) \mathbf 1\notag
\end{align}
Thirdly we have by a rather tedious calculation and  \lemref{binomiallemma}
\begin{align*}
(\text{ad}F_\alpha^{(2)}&(F_\beta K_\beta)\mathbf 1,\text{ad}F_\alpha^{(2)}(F_\beta K_\beta)\mathbf 1)
\\
&=T^4\Big(v^4\left(T^{-2}-T^{-1}[T,-1]+\qbinom{T;-1}{2}\right)[T^2,1]_\beta  \\
 &\quad +v^5(T^{-1}-v[T;-2])[T,-1][T^2;2]_\beta  +v^{6} \qbinom{T;-1}{ 2}[T^2;3]_\beta\Big) .
\end{align*}

Using \corref{firstbetacor}  we get
\begin{align}\label{relation5}
f(-2,-1)^{-2}
&\rho_1\left(r_0\boldsymbol{\beta}^{2,2}_0+r_2\boldsymbol{\beta}^{2,2}_2+r_4\boldsymbol{\beta}^{2,2}_4\right)(\text{ad}F_\alpha(F_\beta^{(2)} K_\beta)\otimes \text{ad}F_\alpha^{(2)}(F_\beta K_\beta))\mathbf 1  \\ 
&=f(2,-1)^{-2}v^4T^{-4}\left(\frac{v^{-2}r_0}{[3]} -
    \frac{r_2v^{-2} }
   {[4]} \qbinom{T;-1}{1}+r_4v^{-4}\left[\begin{matrix} 4
  \\  2 \end{matrix}\right]^{-1} \qbinom{T;-1}{2}\right) \notag
\\ 
&=T^4\Big(v^4\left(T^{-2}-T^{-1}[T,-1]+\qbinom{T;-1}{2}\right)[T^2,1]_\beta  \notag\\
 &\quad +v^5(T^{-1}-v[T;-2])[T,-1][T^2;2]_\beta  +v^{6} \qbinom{T;-1}{ 2}[T^2;3]_\beta\Big)  \notag
\end{align}

From equations \eqnref{relation3}, \eqnref{relation4}  and \eqnref{relation5} we determine the coefficients of $\boldsymbol\beta$:
\begin{align*}
\begin{pmatrix} r_0 \\ r_2[T,-1] \\ r_4 \qbinom{T;-1}{2} \end{pmatrix} =f  \begin{pmatrix}
 v^{-2} & \frac{v^2}{[2]} & v^4 \\
 \frac{[2]}{v^2} & v\{1 \}& -v^2[2]\\
 1 & -v & 1
\end{pmatrix}
\begin{pmatrix}v^4T^2[T^2;1]_\beta   \\ 
	vT^4[2][T^2;1]_\beta +v^6T^7[T;-1]   \\ 
	T^{11}\Big(Tv^5 \{1\}+[T^3,1]_\beta+T^4v^4[T,-1]_\beta\Big) 
\end{pmatrix} \end{align*}
where $f=f(2,-1)^2$.
Here however we have only determined the $\boldsymbol \beta$ for one $U_\nu(\mathfrak a)$-summand of $M(-\varpi_\alpha+\varpi_\beta)$.  In future we plan to determine the coefficients for the other $\boldsymbol\beta$ that appear.

\printindex
\begin{theindex}
\item $\rho_1$,\hfill\pageref{rhoone}
\item $\varrho$,\hfill\pageref{antiauto}
\item $\beta^{m,n}_{2r}$,\hfill \pageref{betamn}
\item $F^{(-k)}\eta$,\hfill \pageref{fminusketa}
\item $U_v(\mathfrak a)$,\hfill\pageref{ua}
\item $\mathcal F_m$,\hfill \pageref{fm}
\item $L$, $L^{-1}$, \hfill\pageref{L}
\item $\mathcal L$,  \hfill\pageref{mathcalLinv}
\item $\mathcal L_{2n}$, \hfill \pageref{mathcalL2n}
\item $\Lambda$, \hfill\pageref{Lambda}
\item $[n]_\alpha$, \hfill\pageref{Gaussintegers}
\item $[T;r]^{(j)},
\qbinom{T;r}{ j}$,\hfill\pageref{quantumbinom}
\end{theindex}
\def\cprime{$'$}

\end{document}